\documentclass{article}
\usepackage{amsmath}
\usepackage{amssymb}
\usepackage{amsthm}
\begin{document}
\def\e#1\e{\begin{equation}#1\end{equation}}
\def\ea#1\ea{\begin{align}#1\end{align}}
\def\eq#1{{\rm(\ref{#1})}}
\theoremstyle{plain}
\newtheorem{thm}{Theorem}[section]
\newtheorem{lem}[thm]{Lemma}
\newtheorem{prop}[thm]{Proposition}
\newtheorem{cor}[thm]{Corollary}
\theoremstyle{definition}
\newtheorem{dfn}[thm]{Definition}
\newtheorem{ex}[thm]{Example}
\def\dim{\mathop{\rm dim}}
\def\Re{\mathop{\rm Re}}
\def\Im{\mathop{\rm Im}}
\def\vol{\mathop{\rm vol}}
\def\U{\mathbin{\rm U}}
\def\SU{\mathop{\rm SU}}
\def\sech{{\textstyle\mathop{\rm sech}}}
\def\ge{\geqslant} 
\def\le{\leqslant} 
\def\R{\mathbin{\mathbb R}}
\def\Z{\mathbin{\mathbb Z}}
\def\C{\mathbin{\mathbb C}}
\def\al{\alpha}
\def\be{\beta}
\def\la{\lambda}
\def\ga{\gamma}
\def\de{\delta}
\def\ep{\epsilon}
\def\th{\theta}
\def\vp{\varphi}
\def\si{\sigma}
\def\ze{\zeta}
\def\Om{\Omega}
\def\Si{\Sigma}
\def\om{\omega}
\def\d{{\rm d}}
\def\pd{\partial}
\def\db{{\bar\partial}}
\def\ts{\textstyle}
\def\sst{\scriptscriptstyle}
\def\w{\wedge}
\def\lt{\ltimes}
\def\sm{\setminus}
\def\bigot{\bigotimes}
\def\iy{\infty}
\def\ra{\rightarrow}
\def\t{\times}
\def\ha{{\textstyle\frac{1}{2}}}
\def\ms#1{\vert#1\vert^2}
\def\md#1{\vert #1 \vert}
\def\bmd#1{\big\vert #1 \big\vert}
\def\cnm#1#2{\Vert #1 \Vert_{C^{#2}}} 
\def\an#1{\langle#1\rangle}
\title{$\U(1)$-invariant special Lagrangian 3-folds. I.
Nonsingular solutions}
\author{Dominic Joyce \\ Lincoln College, Oxford}
\date{}
\maketitle

\section{Introduction}
\label{un1}

Special Lagrangian submanifolds (SL $m$-folds) are a distinguished
class of real $m$-dimensional minimal submanifolds in $\C^m$, which are
calibrated with respect to the $m$-form $\Re(\d z_1\w\cdots\w\d z_m)$.
They can also be defined in (almost) Calabi--Yau manifolds, are
important in String Theory, and are expected to play a r\^ole in the
eventual explanation of Mirror Symmetry between Calabi--Yau 3-folds.

This is the first of a suite of three papers \cite{Joyc6,Joyc7}
studying special Lagrangian 3-folds $N$ in $\C^3$ invariant under
the $\U(1)$-action
\e
{\rm e}^{i\th}:(z_1,z_2,z_3)\mapsto
({\rm e}^{i\th}z_1,{\rm e}^{-i\th}z_2,z_3)
\quad\text{for ${\rm e}^{i\th}\in\U(1)$.}
\label{un1eq1}
\e
These three papers and \cite{Joyc8} are surveyed in \cite{Joyc9}.
Locally we can write $N$ as
\e
\begin{split}
N=\Bigl\{(z_1&,z_2,z_3)\in\C^3:
\Im(z_3)=u\bigl(\Re(z_3),\Im(z_1z_2)\bigr),\\
&\Re(z_1z_2)=v\bigl(\Re(z_3),\Im(z_1z_2)\bigr),
\quad\ms{z_1}-\ms{z_2}=2a\Bigr\},
\end{split}
\label{un1eq2}
\e
where $a\in\R$ and $u,v:\R^2\ra\R$ are differentiable functions.
It will be shown that $N$ is a special Lagrangian 3-fold in $\C^3$
if and only if $u,v$ satisfy
\e
\frac{\pd u}{\pd x}=\frac{\pd v}{\pd y}\quad\text{and}\quad
\frac{\pd v}{\pd x}=-2\bigl(v^2+y^2+a^2\bigr)^{1/2}\frac{\pd u}{\pd y}.
\label{un1eq3}
\e

In fact we have to modify this a bit to allow $N$ to have singularities,
which is one of the main things we are interested in. When $a\ne 0$ it
turns out that $N$ is always nonsingular, and $u,v$ are always smooth
and satisfy \eq{un1eq3} in the usual sense. However, when $a=0$, at
points $(x,0)$ with $v(x,0)=0$ the factor $-2(v^2+y^2+a^2)^{1/2}$ in
\eq{un1eq3} becomes zero, and then \eq{un1eq3} is no longer elliptic.

Because of this, when $a=0$ the appropriate thing to do is to consider
{\it weak solutions} of \eq{un1eq3}, which may have {\it singular points}
$(x,0)$ with $v(x,0)=0$. At such a point $u,v$ may not be differentiable,
and $\bigl(0,0,x+iu(x,0)\bigr)$ is a singular point of the SL 3-fold $N$
in $\C^3$. Weak solutions of \eq{un1eq3} when $a=0$ and their singularities
will be studied in the sequels \cite{Joyc6,Joyc7}, and this paper will
focus on the nonsingular case when~$a\ne 0$.

We begin in \S\ref{un2} with an introduction to special Lagrangian
geometry, and then \S\ref{un3} summarizes some background material
from analysis that we will need later, to do with H\"older spaces
of functions and elliptic operators. Section \ref{un4} considers
special Lagrangian 3-folds invariant under the $\U(1)$-action
\eq{un1eq1}, shows that they can locally be written in the form
\eq{un1eq2} where $u,v$ satisfy \eq{un1eq3}, and gives an explanation
of why \eq{un1eq3} is a {\it nonlinear Cauchy--Riemann equation}
in terms of almost Calabi--Yau geometry. Examples of solutions
$u,v$ of \eq{un1eq3} are given in \S\ref{un5}, and the
corresponding SL 3-folds $N$ in $\C^3$ described.

Section \ref{un6} exploits the fact that \eq{un1eq3} is a nonlinear
Cauchy--Riemann equation, and so $u+iv$ is a bit like a holomorphic
function of $x+iy$. We prove analogues for solutions $u,v$ of \eq{un1eq3}
of well-known results in complex analysis, in particular those involving
multiplicity of zeroes, and formulae counting zeroes of a holomorphic
function in terms of winding numbers.

As an application we show that if $S,T$ are domains in $\R^2$ and
$(\hat u,\hat v):S\ra T$ are solutions of \eq{un1eq3} such that
$\hat u,\hat v,\frac{\pd\hat v}{\pd x}$ and $\frac{\pd\hat v}{\pd y}$
take given values at a point, then there do not exist $(u,v):T\ra
S^\circ$ satisfying \eq{un1eq3} such that $u,v,\frac{\pd v}{\pd x}$
and $\frac{\pd v}{\pd y}$ take given values at a point. This will be
used in \cite{Joyc6} to prove a priori estimates for derivatives of
bounded solutions $u,v$ of \eq{un1eq3} on domains in $\R^2$, and
these in turn will be important in proving the existence of weak
solutions of \eq{un1eq3} when~$a=0$.

In \S\ref{un7} we show that if $S$ is a domain in $\R^2$ and
$u,v\in C^1(S)$ satisfy \eq{un1eq3}, then there exists $f\in C^2(S)$
with $\frac{\pd f}{\pd y}=u$ and $\frac{\pd f}{\pd x}=v$, unique up
to addition of a constant, satisfying
\e
\Bigl(\Bigl(\frac{\pd f}{\pd x}\Bigr)^2+y^2+a^2
\Bigr)^{-1/2}\frac{\pd^2f}{\pd x^2}+2\,\frac{\pd^2f}{\pd y^2}=0.
\label{un1eq4}
\e
This is a {\it second-order quasilinear elliptic equation}. Using
results from analysis, we prove existence and uniqueness of solutions
of the Dirichlet problem for \eq{un1eq4} on strictly convex domains
when $a\ne 0$. Combining this with the results of \S\ref{un4} gives
existence and uniqueness results for nonsingular $\U(1)$-invariant
SL 3-folds in $\C^3$ satisfying certain boundary conditions.

Section \ref{un8} takes a different approach to the same problem.
We show that if $S$ is a domain in $\R^2$ and $u,v\in C^2(S)$ satisfy
\eq{un1eq3}, then $v$ satisfies
\e
\frac{\pd}{\pd x}\Bigl[\bigl(v^2+y^2+a^2\bigr)^{-1/2}
\frac{\pd v}{\pd x}\Bigr]+2\,\frac{\pd^2v}{\pd y^2}=0.
\label{un1eq5}
\e
Again, this is a second-order quasilinear elliptic equation, and we can
prove existence and uniqueness of solutions of the Dirichlet problem for
\eq{un1eq5} on domains in $\R^2$ when $a\ne 0$. This gives existence and
uniqueness results for nonsingular $\U(1)$-invariant SL 3-folds in $\C^3$
satisfying a different kind of boundary condition.

In the sequel \cite{Joyc6} we first prove a priori estimates
for $\frac{\pd u}{\pd x},\frac{\pd u}{\pd y},\frac{\pd v}{\pd x}$ and
$\frac{\pd v}{\pd y}$ when $u,v$ are bounded solutions of \eq{un1eq3}
on a domain $S$ in $\R^2$, and $a\ne 0$. Using these we generalize
Theorems \ref{un7thm2}, \ref{un7thm3}, \ref{un8thm1} and \ref{un8thm2}
below to the case $a=0$, proving existence and uniqueness of {\it weak}\/
solutions $f\in C^1(S)$ and $u,v\in C^0(S)$ to the Dirichlet problems
for \eq{un1eq4} and \eq{un1eq5} on strictly convex domains when $a=0$.
This gives existence and uniqueness results for {\it singular}
$\U(1)$-invariant SL 3-folds in $\C^3$ satisfying certain boundary
conditions.

The following paper \cite{Joyc7} studies these singular solutions
$u,v$ of \eq{un1eq3} when $a=0$ in more detail. We show that under
mild conditions $u,v$ have only isolated singularities, and these
isolated singular points have a {\it multiplicity}, which is a
positive integer, and one of two {\it types}. We also use our results
to construct many {\it special Lagrangian fibrations} on open subsets
of $\C^3$. In \cite{Joyc8} these are used as local models to study
special Lagrangian fibrations of (almost) Calabi--Yau 3-folds, and
to draw some conclusions about the {\it SYZ Conjecture} \cite{SYZ}.
All four papers are reviewed briefly in~\cite{Joyc9}.

A fundamental question about compact special Lagrangian 3-folds $N$
in (almost) Calabi--Yau 3-folds $M$ is: {\it how stable are they under
large deformations\/}? Here we mean both deformations of $N$ in a fixed
$M$, and what happens to $N$ as we deform $M$. The deformation theory
of compact SL 3-folds under {\it small}\/ deformations is already well
understood, and is described in \cite[\S 9]{Joyc4} and \cite[\S 5]{Joyc5}.
But to extend this understanding to large deformations, one needs to take
into account singular behaviour.

One possible moral of this paper and its sequels \cite{Joyc6,Joyc7}
is that {\it compact SL\/ $3$-folds are pretty stable under large
deformations}. That is, we have shown existence and uniqueness for
(possibly singular) $\U(1)$-invariant SL 3-folds in $\C^3$ satisfying
certain boundary conditions. This existence and uniqueness is
{\it entirely unaffected} by singularities that develop in the SL
3-folds, which is quite surprising, as one might have expected that
when singularities develop the existence and uniqueness properties
would break down.

This is encouraging, as both the author's programme for constructing
invariants of almost Calabi--Yau 3-folds in \cite{Joyc1} by counting
special Lagrangian homology 3-spheres, and proving some version of
the SYZ Conjecture \cite{SYZ} in anything other than a fairly weak,
limiting form, will require strong stability properties of compact
SL 3-folds under large deformations; so these papers may be taken
as a small piece of evidence that these two projects may eventually
be successful.
\medskip

\noindent{\it Acknowledgements.} Mark Gross discusses $\U(1)$-invariant
special Lagrangian fibrations in \cite[\S 4]{Gros}, and an idea of his
helped me to make progress at a difficult stage in the composition of
this paper and its sequels. I would also like to thank Rafe Mazzeo and
Rick Schoen for helpful conversations. I was supported by an EPSRC
Advanced Fellowship whilst writing this paper.

\section{Special Lagrangian geometry}
\label{un2}

We now introduce the idea of special Lagrangian submanifolds,
in two different geometric contexts. First, in \S\ref{un21},
we discuss special Lagrangian submanifolds in $\C^m$. Then 
\S\ref{un22} considers special Lagrangian submanifolds in
{\it almost Calabi--Yau manifolds}, K\"ahler manifolds 
equipped with a holomorphic volume form which generalize the
idea of Calabi--Yau manifolds. For an introduction to special
Lagrangian geometry, see Harvey and Lawson \cite{HaLa} or
the author~\cite{Joyc4,Joyc5}.

\subsection{Special Lagrangian submanifolds in $\C^m$}
\label{un21}

We begin by defining {\it calibrations} and {\it calibrated 
submanifolds}, following Harvey and Lawson~\cite{HaLa}.

\begin{dfn} Let $(M,g)$ be a Riemannian manifold. An {\it oriented
tangent\/ $k$-plane} $V$ on $M$ is a vector subspace $V$ of
some tangent space $T_xM$ to $M$ with $\dim V=k$, equipped
with an orientation. If $V$ is an oriented tangent $k$-plane
on $M$ then $g\vert_V$ is a Euclidean metric on $V$, so 
combining $g\vert_V$ with the orientation on $V$ gives a 
natural {\it volume form} $\vol_V$ on $V$, which is a 
$k$-form on~$V$.

Now let $\vp$ be a closed $k$-form on $M$. We say that
$\vp$ is a {\it calibration} on $M$ if for every oriented
$k$-plane $V$ on $M$ we have $\vp\vert_V\le \vol_V$. Here
$\vp\vert_V=\al\cdot\vol_V$ for some $\al\in\R$, and 
$\vp\vert_V\le\vol_V$ if $\al\le 1$. Let $N$ be an 
oriented submanifold of $M$ with dimension $k$. Then 
each tangent space $T_xN$ for $x\in N$ is an oriented
tangent $k$-plane. We say that $N$ is a {\it calibrated 
submanifold} if $\vp\vert_{T_xN}=\vol_{T_xN}$ for all~$x\in N$.
\label{un2def1}
\end{dfn}

It is easy to show that calibrated submanifolds are automatically
{\it minimal submanifolds} \cite[Th.~II.4.2]{HaLa}. Here is the 
definition of special Lagrangian submanifolds in $\C^m$, taken
from~\cite[\S III]{HaLa}.

\begin{dfn} Let $\C^m$ have complex coordinates $(z_1,\dots,z_m)$, 
and define a metric $g$, a real 2-form $\om$ and a complex $m$-form 
$\Om$ on $\C^m$ by
\e
\begin{split}
g=\ms{\d z_1}+\cdots+\ms{\d z_m},\quad
\om&=\frac{i}{2}(\d z_1\w\d\bar z_1+\cdots+\d z_m\w\d\bar z_m),\\
\text{and}\quad\Om&=\d z_1\w\cdots\w\d z_m.
\end{split}
\label{un2eq1}
\e
Then $\Re\Om$ and $\Im\Om$ are real $m$-forms on $\C^m$. Let
$L$ be an oriented real submanifold of $\C^m$ of real dimension 
$m$. We say that $L$ is a {\it special Lagrangian submanifold} 
of $\C^m,$ or {\it SL\/ $m$-fold}\/ for short, if $L$ is calibrated 
with respect to $\Re\Om$, in the sense of Definition~\ref{un2def1}. 
\label{un2def2}
\end{dfn}

As in \cite{Joyc1,Joyc2} there is a more general definition of 
special Lagrangian $m$-fold involving a {\it phase} ${\rm e}^{i\th}$, 
but we will not use it here. Harvey and Lawson \cite[Cor.~III.1.11]{HaLa} 
give the following alternative characterization of special Lagrangian 
submanifolds.

\begin{prop} Let\/ $L$ be a real\/ $m$-dimensional submanifold of\/
$\C^m$. Then $L$ admits an orientation making it into a special
Lagrangian submanifold of\/ $\C^m$ if and only if\/ $\om\vert_L\equiv 0$
and\/~$\Im\Om\vert_L\equiv 0$.
\label{un2prop1}
\end{prop}

An $m$-dimensional submanifold $L$ in $\C^m$ is called {\it Lagrangian} 
if $\om\vert_L\equiv 0$. Thus special Lagrangian submanifolds are 
Lagrangian submanifolds satisfying the extra condition that 
$\Im\Om\vert_L\equiv 0$, which is how they get their name.

Next we give a result characterizing SL 3-planes $\R^3$ in $\C^3$. 
Define an anti-bilinear cross product $\t:\C^3\t\C^3\ra\C^3$ by
\e
(r_1,r_2,r_3)\t(s_1,s_2,s_3)=(\bar r_2\bar s_3-\bar r_3\bar s_2,
\bar r_3\bar s_1-\bar r_1\bar s_3,\bar r_1\bar s_2-\bar r_2\bar s_1).
\label{un2eq2}
\e
It is equivariant under the $\SU(3)$-action on $\C^3$. Using this 
notation, we prove

\begin{prop} Let\/ ${\bf r},{\bf s}\in\C^3$ be linearly independent
over $\R$, with\/ $\om({\bf r},{\bf s})=0$. Then ${\bf r},{\bf s}$
and\/ ${\bf r}\t{\bf s}$ are linearly independent over $\R$, and
$\an{{\bf r},{\bf s},{\bf r}\t{\bf s}}_{\sst\mathbb R}$ is the 
unique special Lagrangian $3$-plane in $\C^3$ 
containing~$\an{{\bf r},{\bf s}}_{\sst\mathbb R}$.
\label{un2prop2}
\end{prop}

\begin{proof} Explicit calculation using \eq{un2eq2} shows that
\begin{gather}
g({\bf r},{\bf r}\t{\bf s})=g({\bf s},{\bf r}\t{\bf s})=0,
\label{un2eq3}\\
\om({\bf r},{\bf r}\t{\bf s})=\om({\bf s},{\bf r}\t{\bf s})=0,
\label{un2eq4}\\
\ms{{\bf r}\t{\bf s}}=\ms{{\bf r}}\ms{{\bf s}}
-g({\bf r},{\bf s})^2-\om({\bf r},{\bf s})^2,
\label{un2eq5}\\
\text{and}\quad (\Im\Om)({\bf r},{\bf s},{\bf r}\t{\bf s})=0,
\label{un2eq6}
\end{gather}
for all ${\bf r},{\bf s}\in\C^3$. When ${\bf r},{\bf s}$ are
linearly independent and $\om({\bf r},{\bf s})=0$, equation
\eq{un2eq3} shows that ${\bf r}\t{\bf s}$ is orthogonal to 
${\bf r},{\bf s}$, and \eq{un2eq5} that $\md{{\bf r}\t{\bf s}}\ne 0$.
Therefore ${\bf r},{\bf s}$ and ${\bf r}\t{\bf s}$ are linearly
independent. 

Also we have $\om({\bf r},{\bf s})=\om({\bf r},{\bf r}\t{\bf s})=
\om({\bf s},{\bf r}\t{\bf s})=0$ by \eq{un2eq4}, so that 
$\an{{\bf r},{\bf s},{\bf r}\t{\bf s}}_{\sst\mathbb R}$ is a 
{\it Lagrangian} 3-plane. Then \eq{un2eq6} shows that
$\an{{\bf r},{\bf s},{\bf r}\t{\bf s}}_{\sst\mathbb R}$ is a 
{\it special}\/ Lagrangian 3-plane, by Proposition \ref{un2prop1}.
It is easy to see that this is the only SL 3-plane in $\C^3$ 
containing~$\an{{\bf r},{\bf s}}_{\sst\mathbb R}$.
\end{proof}

\subsection{Almost Calabi--Yau $m$-folds and SL $m$-folds} 
\label{un22}

We shall define special Lagrangian submanifolds not just in
Calabi--Yau manifolds, as usual, but in the much larger
class of {\it almost Calabi--Yau manifolds}.

\begin{dfn} Let $m\ge 2$. An {\it almost Calabi--Yau $m$-fold}, or
{\it ACY\/ $m$-fold}\/ for short, is a quadruple $(X,J,\om,\Om)$ 
such that $(X,J)$ is a $m$-dimensional complex manifold, $\om$ is 
the K\"ahler form of a K\"ahler metric $g$ on $X$, and $\Om$ is a 
non-vanishing holomorphic $(m,0)$-form on~$X$.

We call $(X,J,\om,\Om)$ a {\it Calabi--Yau $m$-fold}, or {\it CY\/ 
$m$-fold}\/ for short, if in addition $\om$ and $\Om$ satisfy
\e
\om^m/m!=(-1)^{m(m-1)/2}(i/2)^m\Om\w\bar\Om.
\label{un2eq7}
\e
Then for each $x\in X$ there exists an isomorphism $T_xX\cong\C^m$
that identifies $g_x,\om_x$ and $\Om_x$ with the flat versions
$g,\om,\Om$ on $\C^m$ in \eq{un2eq1}. Furthermore, $g$ is Ricci-flat
and its holonomy group is a subgroup of~$\SU(m)$.
\label{un2def3}
\end{dfn}

This is not the usual definition of a Calabi--Yau manifold, but is 
essentially equivalent to it. (Usually one also assumes that $X$ is 
compact). Next, motivated by Proposition \ref{un2prop1}, we define 
special Lagrangian submanifolds of almost Calabi--Yau manifolds.

\begin{dfn} Let $(X,J,\om,\Om)$ be an almost Calabi--Yau $m$-fold 
with metric $g$, and $N$ a real $m$-dimensional submanifold of $X$.
We call $N$ a {\it special Lagrangian submanifold}, or {\it SL\/
$m$-fold\/} for short, if~$\om\vert_N\equiv\Im\Om\vert_N\equiv 0$. 
\label{un2def4}
\end{dfn}

The properties of SL $m$-folds in almost Calabi--Yau $m$-folds
are discussed by the author in \cite{Joyc4,Joyc5}. The deformation and
obstruction theory for {\it compact}\/ SL $m$-folds in almost
Calabi--Yau $m$-folds is well understood, and beautifully behaved.

In this paper we will focus exclusively on special Lagrangian 
3-folds in $\C^3$, and the more general almost Calabi--Yau
context will hardly enter our story at all. However, because
SL $m$-folds in ACY $m$-folds are expected to behave locally 
just like SL $m$-folds in $\C^m$, our results tell us about 
SL 3-folds in ACY 3-folds, especially their singular behaviour. 

\section{Background material from analysis}
\label{un3}

We now briefly summarize some background material we will need
for later analytic results. Our principal reference is Gilbarg 
and Trudinger~\cite{GiTr}.

\subsection{Banach spaces of functions on subsets of $\R^n$}
\label{un31}

We first define a special class of subsets of $\R^n$ called
{\it domains}.

\begin{dfn} A closed, bounded, contractible subset $S$ in $\R^n$ will
be called a {\it domain} if it is a disjoint union $S=S^\circ\cup\pd S$,
where the {\it interior} $S^\circ$ of $S$ is a connected open set in
$\R^n$ with $S=\overline{S^\circ}$, and the {\it boundary} $\pd S=S\sm
S^\circ$ is a compact embedded hypersurface in~$\R^n$.
\label{un3def1}
\end{dfn}

Here the assumption that $S$ is contractible is made for simplicity,
and will not always be necessary. Note that as they are contractible,
domains in $\R^2$ are automatically diffeomorphic to discs. Next we
define some Banach spaces of real functions on~$S$.

\begin{dfn} Let $S$ be a domain in $\R^n$. For each integer $k\ge 0$, 
define $C^k(S)$ to be the space of continuous functions $f:S\ra\R$ with 
$k$ continuous derivatives, and define the norm $\cnm{.}k$ on $C^k(S)$ 
by $\cnm{f}k=\sum_{j=0}^k\sup_S\bmd{\pd^jf}$. Then $C^k(S)$ is a Banach 
space. Define $C^\iy(S)=\bigcap_{k=0}^\iy C^k(S)$ to be the set of smooth 
functions on $S$. It is not a Banach space, with its natural topology.
\label{un3def2}
\end{dfn}

Here $\pd$ is the vector operator $(\frac{\pd}{\pd x_1},\ldots,
\frac{\pd}{\pd x_n})$, where $(x_1,\ldots,x_n)$ are the standard
coordinates on $\R^n$, so that $\pd^jf$ maps $S\ra\bigot^k(\R^n)^*$,
and has components $\frac{\pd^jf}{\pd x_{a_1}\cdots\pd x_{a_j}}$
for $1\le a_1,\ldots,a_j\le n$. The lengths $\bmd{\pd^jf}$ are 
computed using the standard Euclidean metric on~$\R^n$.

\begin{dfn} For $k\ge 0$ an integer and $\al\in(0,1]$, define the 
{\it H\"older space} $C^{k,\al}(S)$ to be the subset of $f\in C^k(S)$ 
for which
\begin{equation*}
[\pd^k f]_\al=\sup_{x\ne y\in S}
\frac{\bmd{\pd^kf(x)-\pd^kf(y)}}{\md{x-y}^\al}
\end{equation*}
is finite, and define the {\it H\"older norm} on $C^{k,\al}(S)$ 
to be $\cnm{f}{k,\al}=\cnm{f}k+[\pd^kf]_\al$. Again, $C^{k,\al}(S)$
is a Banach space.
\label{un3def3}
\end{dfn}

\subsection{Linear and quasilinear elliptic operators}
\label{un32}

We begin by defining {\it second-order linear elliptic operators} 
on functions.

\begin{dfn} Let $S$ be a domain in $\R^n$. A {\it second-order 
linear differential operator} $P$ mapping $C^{k+2}(S)\ra C^k(S)$ or
$C^{k+2,\al}(S)\ra C^{k,\al}(S)$ or $C^\iy(S)\ra C^\iy(S)$
is an operator of the form
\e
\bigl(Pu\bigr)(x)=
\sum_{i,j=1}^na^{ij}(x)\frac{\pd^2u}{\pd x_i\pd x_j}(x)
+\sum_{i=1}^nb^i(x)\frac{\pd u}{\pd x_i}(x)+c(x)u(x),
\label{un3eq1}
\e
where $a^{ij}$, $b^i$ and $c$ lie in $C^k(S)$, or $C^{k,\al}(S)$, or 
$C^\iy(S)$, respectively, and $a^{ij}=a^{ji}$ for all $i,j=1,\ldots,n$.
We call $a^{ij},b^i$ and $c$ the {\it coefficients} of $P$, so 
that, for instance, we say $P$ has $C^{k,\al}$ coefficients 
if $a^{ij}$, $b^i$ and $c$ lie in $C^{k,\al}(S)$. We call $P$ 
{\it elliptic} if the symmetric $n\t n$ matrix $(a^{ij})$ is 
positive definite at every point of~$S$. 
\label{un3def5}
\end{dfn}

There is a much more general definition of ellipticity for 
differential operators of other orders, or acting on vectors
rather than functions, but we will not need it. One can also 
define ellipticity for {\it nonlinear} partial differential 
operators. We will not do this in general, but only for 
{\it quasilinear} differential operators, which are linear 
in their highest-order derivatives. 

\begin{dfn} Let $S$ be a domain in $\R^n$. A {\it second-order 
quasilinear operator} $Q:C^2(S)\ra C^0(S)$ is an operator of the form
\e
\bigl(Qu\bigr)(x)=
\sum_{i,j=1}^na^{ij}(x,u,\pd u)\frac{\pd^2u}{\pd x_i\pd x_j}(x)
+b(x,u,\pd u),
\label{un3eq2}
\e
where $a^{ij}$ and $b$ are continuous maps $S\t\R\t(\R^n)^*\ra\R$,
and $a^{ij}=a^{ji}$ for all $i,j=1,\ldots,n$. We call the functions
$a^{ij}$ and $b$ the {\it coefficients} of $Q$. We call $Q$ 
{\it elliptic} if the symmetric $n\t n$ matrix $(a^{ij})$ is 
positive definite at every point of~$S\t\R\t(\R^n)^*$. 
\label{un3def6}
\end{dfn}

Elliptic operators have good {\it regularity properties} in 
H\"older spaces. 

\begin{thm} Let\/ $S$ be a domain in $\R^n$ and\/ 
$Q:C^2(S)\ra C^0(S)$ a second-order linear or quasilinear 
elliptic differential operator. Suppose that\/ $Qu=f$, with\/
$u\in C^2(S)$ and $f\in C^0(S)$, and\/ $u\vert_{\pd S}=\phi$,
for $\phi\in C^2(\pd S)$. Then
\begin{itemize}
\setlength{\itemsep}{0pt}
\setlength{\parsep}{0pt}
\item[{\rm(a)}] Let\/ $k\ge 0$ and\/ $\al\in(0,1)$, and suppose
that $Q$ has $C^{k,\al}$ coefficients, $f\in C^{k,\al}(S)$, and\/
$\phi\in C^{k+2,\al}(\pd S)$. Then~$u\in C^{k+2,\al}(S)$.
\item[{\rm(b)}] Suppose $Q$ has smooth coefficients, $f\in C^\iy(S)$,
and\/ $\phi\in C^\iy(\pd S)$. Then $u\in C^\iy(S)$.
\item[{\rm(c)}] Suppose $f$ and the coefficients of\/ $Q$ are
real analytic in $S^\circ$. Then $u$ is real analytic in~$S^\circ$.
\end{itemize}
\label{un3thm1}
\end{thm}

\begin{proof} The linear case of part (a) follows from \cite[Th.~6.19, 
p.~111]{GiTr}. For the quasilinear case, regarding $u$ as fixed, write
\begin{equation*}
Pv=\sum_{i,j=1}^na^{ij}(x,u,\pd u)\frac{\pd^2v}{\pd x_i\pd x_j}(x),
\end{equation*}
so that $P$ is a {\it linear} elliptic operator. Applying the linear 
case of (a) to the equation $Pu=f-b(x,u,\pd u)$, we can deduce the 
quasilinear case by induction on $k$. Part (b) follows from (a), 
and part (c) from Morrey~\cite[\S 5.7--\S 5.8]{Morr}.
\end{proof}

Essentially the theorem says that solutions $u$ of an elliptic
equation $Pu=f$ on $S$ are as smooth as possible, given the 
differentiability of $f$ and the boundary condition $\phi$. For
linear elliptic operators $P$ involving only the derivatives of
$u$ there is a {\it maximum principle}~\cite[Th.~3.1, p.~32]{GiTr}:

\begin{thm} Let\/ $S$ be a domain in $\R^n$ and\/ 
$P:C^2(S)\ra C^0(S)$ a second-order linear elliptic 
differential operator of the form \eq{un3eq1}, with
$c(x)\equiv 0$. Suppose $u\in C^0(S)\cap C^2(S^\circ)$.
If\/ $Pu\ge 0$ in $S^\circ$ then the maximum of\/
$u$ is achieved on $\pd S$, and if\/ $Pu\le 0$ in
$S^\circ$ then the minimum of\/ $u$ is achieved on~$\pd S$.
\label{un3thm2}
\end{thm}

\subsection{Existence results for the Dirichlet problem}
\label{un33}

We shall now use results from Gilbarg and Trudinger \cite{GiTr}
to prove existence results for the Dirichlet problem for
two classes of quasilinear elliptic operators, that will be
needed in \S\ref{un7} and \S\ref{un8}. We begin by defining
{\it strictly convex domains} in~$\R^2$.

\begin{dfn} A domain $S$ in $\R^2$ is called {\it strictly 
convex} if $S$ is convex and the curvature of $\pd S$ is
nonzero at every point. So, for example, $x^2+y^2\le 1$ is
strictly convex but $x^4+y^4\le 1$ is not, as its boundary
has zero curvature at $(\pm 1,0)$ and~$(0,\pm 1)$.
\label{un3def7}
\end{dfn}

Here is our first existence result.

\begin{thm} Let\/ $S$ be a strictly convex domain in $\R^2$, and suppose 
\e
\bigl(Pf\bigr)(x)=\sum_{i,j=1}^2a^{ij}(x,f,\pd f)
\frac{\pd^2f}{\pd x_i\pd x_j}(x)
\label{un3eq5}
\e
is a second-order quasilinear elliptic operator in $S$ with\/ 
$a^{ij}\in C^\iy(S\t\R\t\R^2)$. Then whenever $k\ge 0$, $\al\in(0,1)$ 
and\/ $\phi\in C^{k+2,\al}(\pd S)$ there exists a solution $f\in 
C^{k+2,\al}(S)$ of the Dirichlet problem $Pf=0$ in $S$, 
$f\vert_{\pd S}=\phi$. Furthermore $\cnm{f}{1}\le C\cnm{\phi}{2}$ 
for some $C>0$ depending only on~$S$.
\label{un3thm3}
\end{thm}

\begin{proof} It is not difficult to show that as $S$ is strictly
convex there exists $K>0$ depending only on $S$, such that if 
$\phi\in C^2(\pd S)$ then any three distinct points in the graph 
of $\phi$ in $\pd S\t\R\subset\R^2\t\R$ lie in a unique plane in 
$\R^2\t\R$ with slope ${\bf s}\in(\R^2)^*$ satisfying $\md{{\bf s}}
\le K\cnm{\phi}{2}$. In the notation of \cite[p.~310]{GiTr}, the 
boundary data $\pd S,\phi$ satisfies a {\it three point condition}.

Now (noting the equivalence of the three point and bounded slope 
conditions, \cite[p.~314]{GiTr}), \cite[Th.~12.7, p.~312]{GiTr} 
is an existence result for the Dirichlet problem for an operator 
of the form \eq{un3eq5} with boundary data satisfying a three 
point condition. Strengthened as in \cite[Remark (4), p.~314]{GiTr}, 
it implies that if $\phi\in C^{2,\al}(\pd S)$ then there exists 
$f\in C^{2,\al}(S)$ with $Pf=0$ in $S$ and $f\vert_{\pd S}=\phi$,
which satisfies~$\cnm{\pd f}{0}\le K\cnm{\phi}{2}$.

By the maximum principle, Theorem \ref{un3thm2}, the maximum
of $f$ is achieved on $\pd S$. Thus $\cnm{f}{0}=\cnm{\phi}{0}
\le\cnm{\phi}{2}$. Hence
\begin{equation*}
\cnm{f}{1}=\cnm{f}{0}+\cnm{\pd f}{0}\le 
(1+K)\cnm{\phi}{2}=C\cnm{\phi}{2},
\end{equation*}
where $C=1+K$ depends only on $S$. This establishes the case $k=0$ 
of the theorem. If $\phi\in C^{k+2,\al}(S)$ for $k>0$ then $\phi\in 
C^{2,\al}(S)$, so by the $k=0$ case there exists $f\in C^{2,\al}(S)$ 
with $Pf=0$ and $f\vert_{\pd S}=\phi$. But then Theorem \ref{un3thm1} 
shows that $f\in C^{k+2,\al}(S)$, and the proof is complete.
\end{proof}

Combining \cite[Th.~15.12, p.~382]{GiTr} and Theorem \ref{un3thm1} gives:

\begin{thm} Let\/ $S$ be a domain in $\R^n$, and suppose the 
quasilinear operator
\e
\bigl(Qv\bigr)(x)=\sum_{i,j=1}^na^{ij}(x,v)
\frac{\pd^2v}{\pd x_i\pd x_j}(x)+b(x,v,\pd v)
\label{un3eq6}
\e
is elliptic in $S$ with coefficients $a^{ij}\in C^\iy(S\t\R)$ and\/ 
$b\in C^\iy(S\t\R\t\R^n)$ satisfying $\bmd{b(x,v,p)}\le C\ms{p}$ and\/ 
$v\,b(x,v,p)\le 0$ for all\/ $(x,v,p)\in S\t\R\t\R^n$ and some $C>0$.
Then whenever $k\ge 0$, $\al\in(0,1)$ and\/ $\phi\in C^{k+2,\al}(\pd S)$
there exists a solution $v\in C^{k+2,\al}(S)$ of the Dirichlet problem 
$Qv=0$ in $S$, $v\vert_{\pd S}=\phi$.
\label{un3thm4}
\end{thm}

Note that in both theorems, $Q$ is not a general second-order 
quasilinear elliptic operator of the form \eq{un3eq2}, but has some 
restrictions on its structure. In particular, \eq{un3eq5} has $n=2$ 
and no term $b(x,f,\pd f)$, and in \eq{un3eq6} the $a^{ij}$ depend on 
$x$ and $v$ but not on $\pd v$, and the sign of $b$ is restricted.

\section{A class of $\U(1)$-invariant SL 3-folds in $\C^3$}
\label{un4}

We will now study special Lagrangian 3-folds $N$ in $\C^3$ 
invariant under the $\U(1)$-action
\e
{\rm e}^{i\th}:(z_1,z_2,z_3)\mapsto
({\rm e}^{i\th}z_1,{\rm e}^{-i\th}z_2,z_3)
\quad\text{for ${\rm e}^{i\th}\in\U(1)$.}
\label{un4eq1}
\e

We shall assume that $N$ may be written
\e
\begin{split}
N=\Bigl\{(z_1&,z_2,z_3)\in\C^3:
\Im(z_3)=u\bigl(\Re(z_3),\Im(z_1z_2)\bigr),\\
&\Re(z_1z_2)=v\bigl(\Re(z_3),\Im(z_1z_2)\bigr),
\quad\ms{z_1}-\ms{z_2}=2a\Bigr\},
\end{split}
\label{un4eq2}
\e
where $a\in\R$ and $u,v:\R^2\ra\R$ are continuous functions, which are
smooth except perhaps at certain singular points. Here is why we choose
to write $N$ in this form. As the functions $\Re(z_1z_2),\Im(z_1z_2),
\ms{z_1}-\ms{z_2},\Re(z_3)$ and $\Im(z_3)$ involved in \eq{un4eq2} 
are $\U(1)$-invariant, $N$ is automatically $\U(1)$-invariant. 

Also, as in \cite[Prop.~4.2]{Joyc2}, if $N$ is a connected Lagrangian 
submanifold of $\C^m$ invariant under a Lie subgroup $G$ of the 
automorphism group $\U(m)\lt\C^m$ of $\C^m$, then the moment map $\mu$ 
of $G$ is constant on $N$. Now the moment map of the $\U(1)$-action 
\eq{un4eq1} is $\ms{z_1}-\ms{z_2}$. Thus $\ms{z_1}-\ms{z_2}=2a$ for
some $a\in\R$ on any $\U(1)$-invariant SL 3-fold $N$ in $\C^3$, which
is why we have taken $\ms{z_1}-\ms{z_2}=2a$ to be one of the equations
defining~$N$.

In the other two equations $\Re(z_1z_2)=v\bigl(\Re(z_3),\Im(z_1z_2)\bigr)$
and $\Im(z_3)=u\bigl(\Re(z_3),\Im(z_1z_2)\bigr)$, what we are doing is
regarding the functions $x=\Re(z_3)$ and $y=\Im(z_1z_2)$ as 
{\it coordinates} on $N/\U(1)$, and expressing the other two degrees of 
freedom $\Re(z_1z_2)$ and $\Im(z_3)$ as functions of $x$ and $y$. Thus
we define $N$ as a kind of graph of the pair of functions~$(u,v)$.

Note that not every $\U(1)$-invariant SL 3-fold $N$ in $\C^3$ may be 
written in the form \eq{un4eq2}. Locally this is generally possible,
but globally the functions $u$ and $v$ would have to be multi-valued,
branched covers of $\R^2$ for instance. However, we will see that the
class of SL 3-folds of this form do have many nice properties, and are
interesting both in themselves and for our later applications. So 
equation \eq{un4eq2} should be regarded as more than just an arbitrary 
choice of coordinate system.

\subsection{Finding the equations on $u$ and $v$}
\label{un41}

We now calculate the conditions on the functions $u(x,y)$, $v(x,y)$
for the 3-fold $N$ of \eq{un4eq2} to be special Lagrangian.

\begin{prop} Let\/ $S$ be a domain in $\R^2$ or $S=\R^2$, let\/
$u,v:S\ra\R$ be continuous, and\/ $a\in\R$. Define
\e
\begin{split}
N=\bigl\{(z_1,z_2,z_3)\in\C^3:\,& z_1z_2=v(x,y)+iy,\quad z_3=x+iu(x,y),\\
&\ms{z_1}-\ms{z_2}=2a,\quad (x,y)\in S\bigr\}.
\end{split}
\label{un4eq3}
\e
Then 
\begin{itemize}
\setlength{\itemsep}{0pt}
\setlength{\parsep}{0pt}
\item[{\rm(a)}] If\/ $a=0$, then $N$ is a (possibly singular)
special Lagrangian $3$-fold in $\C^3$, with boundary over
$\pd S$, if\/ $u,v$ are differentiable and satisfy
\e
\frac{\pd u}{\pd x}=\frac{\pd v}{\pd y}
\quad\text{and}\quad
\frac{\pd v}{\pd x}=-2\bigl(v^2+y^2\bigr)^{1/2}\frac{\pd u}{\pd y},
\label{un4eq4}
\e
except at points $(x,0)$ in $S$ with\/ $v(x,0)=0$, where $u,v$ 
need not be differentiable. The singular points of\/ $N$ are those
of the form $(0,0,z_3)$, where $z_3=x+iu(x,0)$ for $x\in\R$ 
with\/~$v(x,0)=0$.
\item[{\rm(b)}] If\/ $a\ne 0$, then $N$ is a nonsingular SL\/
$3$-fold in $\C^3$, with boundary over $\pd S$, if and only if\/
$u,v$ are differentiable on all of\/ $S$ and satisfy
\e
\frac{\pd u}{\pd x}=\frac{\pd v}{\pd y}\quad\text{and}\quad
\frac{\pd v}{\pd x}=-2\bigl(v^2+y^2+a^2\bigr)^{1/2}\frac{\pd u}{\pd y}.
\label{un4eq5}
\e
\end{itemize}
\label{un4prop1}
\end{prop}

\begin{proof} We shall give the proof for part (a). Part (b) is similar 
but more complicated, and will be left to the reader. Let $a=0$, let $N$ 
be defined by \eq{un4eq3}, and let ${\bf z}=(z_1,z_2,z_3)\in N$. For
$\bf z$ to be a nonsingular point of $N$, we need $u$ and $v$ to be
differentiable at $(x,y)=\bigl(\Re(z_3),\Im(z_1z_2)\bigr)$ in $S$,
and for the derivatives of the three functions
\begin{equation*}
\Re(z_1z_2)-v\bigl(\Re(z_3),\Im(z_1z_2)\bigr),\;\>
\Im(z_3)-u\bigl(\Re(z_3),\Im(z_1z_2)\bigr),\;\>\ms{z_1}-\ms{z_2}
\end{equation*}
on $\C^3$ to be linearly independent at~$\bf z$.

Now if $z_1=z_2=0$ then $\ms{z_1}-\ms{z_2}$ has zero derivative at 
$\bf z$. Thus points of the form $(0,0,z_3)$ in $N$ will be singular. 
Clearly, these occur exactly when $z_3=x+iu(x,0)$ for $x\in\R$ with 
$v(x,0)=0$. Also, as $\ms{z_1}-\ms{z_2}=0$, such points occur in $N$ 
only when $a=0$. We shall see that these are the only singular points
in $N$, provided $u$ and $v$ are differentiable. 

To prove part (a) we need to show that each ${\bf z}\in N$ not
of the form $(0,0,z_3)$ is a nonsingular point of $N$, and the
tangent space $T_{\bf z}N$ is a special Lagrangian 3-plane $\R^3$ 
in $\C^3$. As $N$ is $\U(1)$-invariant, it is enough to prove this
for one point in each orbit of the $\U(1)$-action \eq{un4eq1}.
Since $\md{z_1}=\md{z_2}$ on $N$, each $\U(1)$-orbit in $N$ contains 
one or two points $(z_1,z_2,z_3)$ with~$z_1=z_2$. 

Thus it is enough to show that $T_{\bf z}N$ exists and is special 
Lagrangian for points ${\bf z}=(z_1,z_1,z_3)$ in $N$ with $z_1\ne 0$. 
In our next lemma we identify $T_{\bf z}N$ at such a point. The proof 
is elementary, and is left as an exercise.

\begin{lem} Let\/ ${\bf z}=(z_1,z_1,z_3)\in N$, with\/ $z_1\ne 0$.
Set\/ $x=\Re(z_3)$ and\/ $y=\Im(z_1^2)$. Then $N$ is nonsingular
at\/ $\bf z$, and\/ $T_{\bf z}N=\an{{\bf p}_1,{\bf p}_2,
{\bf p}_3}_{\sst\mathbb R},$ where
\ea
{\bf p}_1&=(iz_1,-iz_1,0),
\label{un4eq6}\\
{\bf p}_2&=\bigl((2z_1)^{-1}{\ts\frac{\pd v}{\pd x}}(x,y),
(2z_1)^{-1}{\ts\frac{\pd v}{\pd x}}(x,y),
1+i{\ts\frac{\pd u}{\pd x}}(x,y)\bigr) \quad\text{and}
\label{un4eq7}\\
{\bf p}_3&=\bigl((2z_1)^{-1}({\ts\frac{\pd v}{\pd y}}(x,y)+i),
(2z_1)^{-1}({\ts\frac{\pd v}{\pd y}}(x,y)+i),
i{\ts\frac{\pd u}{\pd y}}(x,y)\bigr).
\label{un4eq8}
\ea
\label{un4lem1}
\end{lem}

Now define $\t:\C^3\t\C^3\ra\C^3$ as in \eq{un2eq2}, and 
apply Proposition \ref{un2prop2} with ${\bf r}={\bf p}_1$ and 
${\bf s}={\bf p}_2$. Clearly ${\bf p}_1$ and ${\bf p}_2$ are 
linearly independent, and $\om({\bf p}_1,{\bf p}_2)=0$. So
Proposition \ref{un2prop2} shows that $\an{{\bf p}_1,{\bf p}_2,
{\bf p}_1\t{\bf p}_2}_{\sst\mathbb R}$ is the unique SL 3-plane
in $\C^3$ containing~$\an{{\bf p}_1,{\bf p}_2}_{\sst\mathbb R}$. 

Therefore $\an{{\bf p}_1,{\bf p}_2,{\bf p}_3}_{\sst\mathbb R}$ is
an SL 3-plane if and only if ${\bf p}_3\in\an{{\bf p}_1,{\bf p}_2,
{\bf p}_1\t{\bf p}_2}_{\sst\mathbb R}$. Combining equations
\eq{un2eq2}, \eq{un4eq6} and \eq{un4eq7} gives
\e
{\bf p}_1\t{\bf p}_2=\bigl(\bar z_1({\ts\frac{\pd u}{\pd x}}+i),
\bar z_1({\ts\frac{\pd u}{\pd x}}+i),-i{\ts\frac{\pd v}{\pd x}}\bigr).
\label{un4eq9}
\e
So suppose ${\bf p}_3=\al{\bf p}_1+\be{\bf p}_2+\ga{\bf p}_1\t{\bf p}_2$.
As the first two coordinates are equal in ${\bf p}_2,{\bf p}_3$ and
${\bf p}_1\t{\bf p}_2$ but not in ${\bf p}_1$, we see that $\al=0$.
Taking real parts in the third coordinate gives $\be=0$. And comparing
real multiples of $i\bar z_1$ in the first coordinate shows 
that~$\ga=\ha\md{z_1}^{-2}$.

Thus $T_{\bf z}N$ is special Lagrangian if and only if
${\bf p}_1\t{\bf p}_2=2\ms{z_1}{\bf p}_3$. By \eq{un4eq8} and
\eq{un4eq9}, this reduces to 
\e
\frac{\pd u}{\pd x}=\frac{\pd v}{\pd y}\quad\text{and}\quad
\frac{\pd v}{\pd x}=-2\ms{z_1}\frac{\pd u}{\pd y}\quad\text{at $(x,y)$.}
\label{un4eq10}
\e
But $v=\Re(z_1^2)$ and $y=\Im(z_1^2)$ by \eq{un4eq3}, so that 
$\md{z_1}^4=v^2+y^2$, and $\ms{z_1}=(v^2+y^2)^{1/2}$. Substituting 
this into \eq{un4eq10} gives equation \eq{un4eq4}, which proves part 
(a) of Proposition \ref{un4prop1}. Part (b) is left to the reader.
\end{proof}

Equations \eq{un4eq4} and \eq{un4eq5} are {\it nonlinear versions
of the Cauchy--Riemann equations}. For if we replace the factors 
$2(v^2+y^2)^{1/2}$ and $2(v^2+y^2+a^2)^{1/2}$ in \eq{un4eq4} and 
\eq{un4eq5} by 1, the equations become
\begin{equation*}
\frac{\pd u}{\pd x}=\frac{\pd v}{\pd y}\quad\text{and}\quad
\frac{\pd v}{\pd x}=-\,\frac{\pd u}{\pd y},
\end{equation*}
which are the conditions for $u+iv$ to be a holomorphic function of 
$x+iy$. We may therefore expect the solutions of \eq{un4eq4} and 
\eq{un4eq5} to have qualitative features in common with solutions of 
the Cauchy--Riemann equations.

\begin{prop} Let\/ $S$ be a domain in $\R^2$, let\/ $a\ne 0$, and 
suppose $u,v\in C^1(S)$ satisfy \eq{un4eq5}. Then $u,v$ are real 
analytic in $S^\circ$, and satisfy
\ea
\frac{\pd^2u}{\pd x^2}+2\bigl(v^2\!+\!y^2\!+\!a^2\bigr)^{1/2}
\frac{\pd^2 u}{\pd y^2}+2\,\frac{\pd}{\pd y}\Bigl[
\bigl(v^2\!+\!y^2\!+\!a^2\bigr)^{1/2}\Bigr]\frac{\pd u}{\pd y}&=0
\quad\text{and}
\label{un4eq11}\\
\bigl(v^2\!+\!y^2\!+\!a^2\bigr)^{-1/2}\frac{\pd^2 v}{\pd x^2}
+2\,\frac{\pd^2v}{\pd y^2}+\frac{\pd}{\pd x}\Bigl[
\bigl(v^2\!+\!y^2\!+\!a^2\bigr)^{-1/2}\Bigr]\frac{\pd v}{\pd x}&=0
\quad\text{in $S^\circ$.}
\label{un4eq12}
\ea
\label{un4prop3}
\end{prop}

\begin{proof} One can show that $u,v$ are real analytic in $S^\circ$
following Harvey and Lawson \cite[Th.~III.2.7]{HaLa}. Thus $v$ is 
twice continuously differentiable, so that $\frac{\pd}{\pd x}\bigl[
\frac{\pd v}{\pd y}\bigr]=\frac{\pd}{\pd y}\bigl[\frac{\pd v}{\pd x}
\bigr]$ in $S^\circ$. Using \eq{un4eq5} to substitute for 
$\frac{\pd v}{\pd y},\frac{\pd v}{\pd x}$ in terms of 
$\frac{\pd u}{\pd x},\frac{\pd u}{\pd y}$ gives \eq{un4eq11}. 
Equation \eq{un4eq12} follows in the same way.
\end{proof}

Regarding the factors $(v^2+y^2+a^2)^{\pm 1/2}$ as part of the
coefficients $a^{ij}(x),b^i(x)$, we see that \eq{un4eq11} and 
\eq{un4eq12} are second-order linear elliptic equations in $u$ 
and $v$ respectively, of the form \eq{un3eq1}, with $c(x)\equiv 0$. 
Therefore by the maximum principle, Theorem \ref{un3thm2}, we have:

\begin{cor} Let\/ $S$ be a domain in $\R^2$, let\/ $a\ne 0$,
and suppose $u,v\in C^1(S)$ satisfy \eq{un4eq5}. Then the
maxima and minima of\/ $u$ and\/ $v$ are achieved on~$\pd S$.
\label{un4cor}
\end{cor}

\subsection{Interpretation using K\"ahler quotients}
\label{un42}

We can use an idea due independently to Goldstein \cite[\S 2]{Gold} 
and Gross \cite[\S 1]{Gros} to interpret some features of the 
above construction. Let $(X,J,\om,\Om)$ be an almost Calabi--Yau
$m$-fold, as in \S\ref{un22}, and $G$ a $k$-dimensional Lie group 
acting on $X$ preserving $J,\om,\Om$, with Lie algebra $\mathfrak g$. 
Suppose the $G$-action admits a moment map~$\mu:X\ra{\mathfrak g}^*$. 

Then for each $c\in Z({\mathfrak g}^*)$, the quotient $M_c=\mu^{-1}(c)/G$ 
is nonsingular wherever $G$ acts freely, and has the structure of an 
almost Calabi--Yau $(m\!-\!k)$-fold on its nonsingular part. If $N$ is 
a connected, $G$-invariant SL $m$-fold in $X$, then $N\subset\mu^{-1}(c)$ 
for some $c\in Z({\mathfrak g}^*)$, and $L=N/G$ is an SL $(m\!-\!k)$-fold 
in $M_c$. Conversely, if $L$ is an SL $(m\!-\!k)$-fold in $M_c$ then $L$ 
pulls back to an SL $m$-fold $N$ in $X$, contained in~$\mu^{-1}(c)$.

In our case, $X$ is $\C^3$ and $G$ is $\U(1)$, acting as in 
\eq{un4eq1}. Any $\U(1)$-invariant SL 3-fold $N$ in $\C^3$
lies in $\mu^{-1}(2a)$ for some $a\in\R$, where $\mu(z_1,z_2,z_3)
=\ms{z_1}-\ms{z_2}$, and pushes down to an SL 2-fold 
in~$M_a=\mu^{-1}(2a)/\U(1)$.

Now SL 2-folds in an almost Calabi--Yau 2-fold $(M,I,\om,\Om)$ are 
the same thing as {\it pseudoholomorphic curves} in $M$ with respect 
to an alternative almost complex structure $J$ depending on $I,\om$ 
and $\Om$. Thus, finding $\U(1)$-invariant SL 3-folds $N$ in $\C^3$ 
is equivalent to finding pseudoholomorphic curves $\Si$ in a family 
of almost complex 2-folds~$M_a$. 

Therefore, it is not surprising that \eq{un4eq4} and \eq{un4eq5} are
nonlinear versions of the Cauchy--Riemann equations. However, this
almost complex point of view is not that helpful in understanding
the {\it singular points} of $N$, which occur when $a=v=y=0$. For the
$\U(1)$-action on $\mu^{-1}(0)$ is not free, and thus $M_0=\mu^{-1}(0)/
\U(1)$ is a {\it singular} almost complex 2-fold.

So the problem is not one of studying singular pseudoholomorphic 
curves in a nonsingular almost complex 2-fold, which are already 
very well understood, but of studying pseudoholomorphic curves in 
a singular almost complex 2-fold, where the almost complex structure 
itself has unpleasant, non-isolated singularities, which are not at
all like the singularities of complex manifolds.

\section{Examples}
\label{un5}

By starting with known examples $N$ of SL 3-folds in $\C^3$
invariant under the $\U(1)$-action \eq{un4eq1} and solving
\eq{un4eq2} for $u$ and $v$, we can construct examples of
solutions $u,v$ to equations \eq{un4eq4} and~\eq{un4eq5}.

We shall do this with a family of explicit SL 3-folds in $\C^3$ 
written down by Harvey and Lawson \cite[\S III.3.A]{HaLa}, and
studied in more detail by the author \cite[\S 3]{Joyc1}. Let 
$a\ge 0$. Define a subset $N_a$ in $\C^3$ by
\e
\begin{split}
N_a=\Bigl\{(z_1&,z_2,z_3)\in\C^3:\ms{z_1}-2a=\ms{z_2}=\ms{z_3},\\
&\Im\bigl(z_1z_2z_3\bigr)=0,\quad \Re\bigl(z_1z_2z_3\bigr)
\ge 0\Bigr\}.
\end{split}
\label{un5eq1}
\e
By \cite[\S III.3.A]{HaLa} and \cite[\S 3]{Joyc1}, $N_a$ is a
nonsingular SL 3-fold diffeomorphic to ${\mathcal S}^1\t\R^2$ when 
$a>0$, and $N_0$ is an SL $T^2$-cone with one singular point at
$(0,0,0)$. We shall show that these SL 3-folds can be written in 
the form~\eq{un4eq2}.

\begin{thm} Let\/ $a\ge 0$. Then there exist unique
$u_a,v_a:\R^2\ra\R$ such that
\e
\begin{split}
N=\Bigl\{(z_1&,z_2,z_3)\in\C^3:
\Im(z_3)=u_a\bigl(\Re(z_3),\Im(z_1z_2)\bigr),\\
&\Re(z_1z_2)=v_a\bigl(\Re(z_3),\Im(z_1z_2)\bigr),
\quad\ms{z_1}-\ms{z_2}=2a\Bigr\}
\end{split}
\label{un5eq2}
\e
is the special Lagrangian $3$-fold\/ $N_a$ of\/ \eq{un5eq1}. Furthermore:
\begin{itemize}
\setlength{\itemsep}{0pt}
\setlength{\parsep}{0pt}
\item[{\rm(a)}] $u_a,v_a$ are smooth on $\R^2$ and satisfy \eq{un4eq5}, 
except at\/ $(0,0)$ when $a=0$, where they are only continuous.
\item[{\rm(b)}] $u_a(x,y)<0$ when $y>0$ for all\/ $x$, and\/ $u_a(x,0)=0$
for all\/ $x$, and\/ $u_a(x,y)>0$ when $y<0$ for all\/~$x$.
\item[{\rm(c)}] $v_a(x,y)>0$ when $x>0$ for all\/ $y$, and\/ $v_a(0,y)=0$
for all\/ $y$, and\/ $v_a(x,y)<0$ when $x<0$ for all\/~$y$.
\item[{\rm(d)}] $u_a(0,y)=-y\bigl(\md{a}+\sqrt{y^2+a^2}\,\,\bigr)^{-1/2}$ 
for all\/~$y$.
\item[{\rm(e)}] $v_a(x,0)=x\bigl(x^2+2\md{a}\bigr)^{1/2}$ for all\/~$x$.
\end{itemize}
\label{un5thm}
\end{thm}

\begin{proof} For simplicity, we first consider the case $a=0$. Let 
$N_0$ be as in \eq{un5eq1}, let $(z_1,z_2,z_3)\in N_0$, and set
\e
x=\Re(z_3),\quad y=\Im(z_1z_2),\quad u=\Im z_3
\quad\text{and}\quad v=\Re(z_1z_2).
\label{un5eq3}
\e
Then $z_3=x+iu$, and $z_1z_2=v+iy$. Thus the first condition 
$\ms{z_1}=\ms{z_2}=\ms{z_3}$ in \eq{un5eq1} becomes
\begin{equation*}
\ms{z_1}=\ms{z_2}=x^2+u^2.
\end{equation*}
Squaring gives $\ms{z_1z_2}=(x^2+u^2)^2$, so 
substituting for $z_1z_2$ yields
\e
v^2+y^2=(x^2+u^2)^2.
\label{un5eq4}
\e
Similarly, using the expressions for $z_1z_2$ and $z_3$ above, the 
second and third conditions on $(z_1,z_2,z_3)$ in \eq{un5eq1} become
\e
vu+yx=0 \qquad\text{and}\qquad vx-yu\ge 0.
\label{un5eq5}
\e

We will use equations \eq{un5eq4} and \eq{un5eq5} to prove
parts (b) and (c) of the theorem. First suppose $y=0$.
Then \eq{un5eq5} gives $vu=0$, so $v=0$ or $u=0$. If $v=0$ 
then \eq{un5eq4} gives $x^2+u^2=0$, so $x=u=0$. Thus $y=0$ 
implies $u=0$. Similarly $u=0$ implies $y=0$, so $u=0$ if 
and only if $y=0$, as in part (b). In the same way $v=0$ 
if and only if $x=0$, as in part~(c). 

We claim that the two terms $vx$ and $-yu$ in \eq{un5eq5} 
are both nonnegative. If one is zero this is obvious. So suppose 
both are nonzero, so that $x,y,u$ and $v$ are all nonzero. From 
\eq{un5eq5}, the signs of three of these terms determine the sign of 
the fourth. It is easy to verify that for all eight sign possibilities,
$vx$ and $-yu$ have the same sign. So both are nonnegative
by \eq{un5eq5}. Hence $yu\le 0$, and $u=0$ if and only if $y=0$. 
Clearly, this proves part (b). Part (c) follows in the same way.

Next we shall show that for each pair $(x,y)$, there is exactly 
one pair $(u,v)$ satisfying \eq{un5eq4} and \eq{un5eq5}. Multiplying
\eq{un5eq4} by $u^2$ and replacing $v^2u^2$ by $y^2x^2$ using 
\eq{un5eq5}, we get $u^6+2x^2u^4+(x^2-y^2)u^2-y^2x^2=0$. This is a 
sextic in $u$, independent of $v$. Putting $\al=u^2$, it becomes 
\begin{equation*}
P(\al)=\al^3+2x^2\al^2+(x^2-y^2)\al-y^2x^2=0.
\end{equation*}

Thus $u^2$ is a real, nonnegative root of the cubic $P$. Divide
into cases
\begin{itemize}
\setlength{\itemsep}{0pt}
\setlength{\parsep}{0pt}
\item[(i)] $x\ne 0$, $y\ne 0$ and $P$ has three real roots 
$\ga_1,\ga_2,\ga_3$, not necessarily distinct;
\item[(ii)] $x\ne 0$, $y\ne 0$ and $P$ has one real root
$\ga$ and a complex conjugate pair of non-real roots~$\de,\bar\de$;
\item[(iii)] $y=0$; and (iv) $x=0$ and $y\ne 0$.
\end{itemize}
We shall show that in cases (i)--(iii), the cubic $P$ has exactly 
one real nonnegative root, giving a unique value of $u^2$. In case
(iv) there are two nonnegative roots, but one can be excluded.

In case (i) we have $\ga_1+\ga_2+\ga_3=-2x^2<0$, so at least
one $\ga_j$ is negative. But $\ga_1\ga_2\ga_3=y^2x^2>0$, so an
even number of $\ga_j$ are negative and an odd number positive. The
only possibility is that one $\ga_j$ is positive and two negative.
So $P$ has exactly one nonnegative root. In case (ii) we have 
$\ga\ms{\de}=y^2x^2>0$, proving that $\ga>0$, so $P$ has exactly 
one nonnegative root. In case (iii) we have $P(\al)=\al\bigl(\al+
x^2\bigr)^2$, with roots 0 and $-x^2$ (twice), so the only 
nonnegative root is~0.

In case (iv) we have $P(\al)=\al^3-y^2\al$, with roots $y,0$ and 
$-y$. Thus there are two nonnegative roots, $\md{y}$ and 0. However, 
if $\al=0$ then $u^2=0$, and $x^2=0$ by assumption, so the
right hand side of \eq{un5eq4} is zero. But $y\ne 0$, so the left
hand side is positive, a contradiction. Hence $\al\ne 0$, and there
is one allowable value for $\al$, which is~$\md{y}$.

We have shown that \eq{un5eq4} and \eq{un5eq5} determine $u^2$ 
uniquely, and that there is a solution $u^2$ for all $x,y$. 
This yields $u$ up to sign. But part (b) gives the sign of 
$u$, so $u$ is determined uniquely. If $u\ne 0$, equation 
\eq{un5eq5} determines $v$. If $u=0$ then $y=0$ by (b), so
\eq{un5eq4} gives $v^2=x^2$, and $v=\pm x$. The sign
of $v$ is given by (c). Therefore for all pairs $x,y$, there
are unique solutions $u,v$ to \eq{un5eq4} and~\eq{un5eq5}.

Let us review what we have proved so far. If $(z_1,z_2,z_3)\in
N_0$ and $x,y,u,v$ are defined by \eq{un5eq3}, then they
satisfy \eq{un5eq4} and \eq{un5eq5}. Also, given any $x,y$ there
exist unique $u,v$ satisfying \eq{un5eq4} and \eq{un5eq5}. So, 
putting $u_0(x,y)=u$ and $v_0(x,y)=v$ defines the functions 
$u_0,v_0$ in the theorem uniquely, and then $N_0$ is a subset 
of the 3-fold $N$ of \eq{un5eq2}. The converse, that $N\subseteq N_0$, 
follows easily by reversing the argument above, since if $(z_1,z_2,z_3)
\in N$ then \eq{un5eq4} and \eq{un5eq5} are equivalent to the equations 
defining $N_0$. Hence~$N=N_0$. 

It remains to prove parts (a), (d) and (e). The smoothness in (a) 
follows directly from \eq{un5eq4} and \eq{un5eq5}, or indirectly from 
the fact that $N_0$ is smooth except at $(0,0,0)$, and $u_0,v_0$ satisfy 
\eq{un4eq5} where they are smooth by Proposition \ref{un4prop1}. 
For part (d), set $x=0$. Then $v=0$ by (c), so \eq{un5eq4} gives 
$u^4=y^2$. So $u_0(0,y)=\pm\md{y}^{1/2}$, and the sign is determined 
by (b). Part (e) follows in the same way. 
This completes the proof for~$a=0$. 

When $a\ne 0$, equation \eq{un5eq4} must be replaced by
\begin{equation*}
v^2+y^2=\bigl(x^2+u^2\bigr)\bigl(x^2+u^2+2\md{a}\bigr),
\end{equation*}
but the rest of the proof is more-or-less unchanged.
\end{proof}

Here are some remarks on the theorem.
\begin{itemize}
\setlength{\itemsep}{0pt}
\setlength{\parsep}{0pt}
\item Let $a>0$. As \eq{un4eq5} depends only on $a^2$, the functions 
$u_a,v_a$ also solve \eq{un4eq5} with $a$ replaced by $-a$. The 
corresponding SL 3-fold is 
\begin{align*}
N_{-a}=\Bigl\{(z_1,z_2,z_3)\in\C^3:\,&\ms{z_1}=\ms{z_2}-2a=\ms{z_3},\\
&\Im\bigl(z_1z_2z_3\bigr)=0,\quad \Re\bigl(z_1z_2z_3\bigr)\ge 0\Bigr\}.
\end{align*}
\item The SL 3-fold $N_0$ is a {\it cone} in $\C^3$, so that 
$tN_0=N_0$ for all $t>0$. It follows that the functions 
$u_0,v_0$ constructed above satisfy
\e
u_0(tx,t^2y)=tu_0(x,y)\;\>\text{and}\;\>
v_0(tx,t^2y)=t^2v_0(x,y) \;\>\text{for all $t>0$,}
\label{un5eq6}
\e
a kind of {\it weighted homogeneity equation}. 
\item The functions $u_0,v_0$ in the theorem are not smooth at 
$(0,0)$. Their behaviour helps us to guess properties of more 
general singular solutions to \eq{un4eq4}. For instance,
$u_0(0,y)=y\md{y}^{-1/2}$ by (d), so $\frac{\pd u_0}{\pd y}$ is 
unbounded near $(0,0)$. This will be important when we consider 
the problem of finding {\it a priori estimates} for derivatives
of solutions $u,v$ of \eq{un4eq4} in~\cite{Joyc6}.
\end{itemize}

Here are some other explicit examples of solutions to \eq{un4eq4} 
and~\eq{un4eq5}. 

\begin{ex} Let $\al,\be,\ga\in\R$ and define $u(x,y)=\al x+\be$ 
and $v(x,y)=\al y+\ga$. Then $u,v$ satisfy \eq{un4eq5} for any 
value of~$a$. 
\label{un5ex1}
\end{ex}

\begin{ex} Let $S=\R^2$, $u(x,y)=y\tanh x$ and $v(x,y)=\ha
y^2\sech^2x-\ha\cosh^2x$. Then $u$ and $v$ satisfy \eq{un4eq4}.
Equation \eq{un4eq3} with $a=0$ defines an explicit nonsingular
special Lagrangian 3-fold $N$ in $\C^3$. It can be shown that $N$
is ruled, and arises from Harvey and Lawson's `austere submanifold'
construction \cite[\S III.3.C]{HaLa} of SL $m$-folds in $\C^m$, as
the normal bundle of a catenoid in~$\R^3$.
\label{un5ex2}
\end{ex}

\begin{ex} Let $S=\R^2$, $u(x,y)=\md{y}-\ha\cosh 2x$ and
$v(x,y)=-y\sinh 2x$. Then $u,v$ satisfy \eq{un4eq4}, except that
$\frac{\pd u}{\pd y}$ is not well-defined on the line $y=0$. So
equation \eq{un4eq3} defines an explicit special Lagrangian 3-fold
$N$ in $\C^3$. It turns out that $N$ is the union of two nonsingular
SL 3-folds intersecting in a real curve, which are constructed in
\cite[Ex.~7.4]{Joyc3} by evolving paraboloids in~$\C^3$.
\label{un5ex3}
\end{ex}

\section{Results using `winding number' techniques}
\label{un6}

We will now discuss some results based on the idea of
{\it winding number}.

\begin{dfn} Let $C$ be a compact oriented 1-manifold, and
$\ga:C\ra\R^2\sm\{0\}$ a differentiable map. Then the
{\it winding number of\/ $\ga$ about\/ $0$ along} $C$ is
$\frac{1}{2\pi}\int_C\ga^*(\d\th)$, where $\d\th$ is the
closed 1-form $(x\,\d y-y\,\d x)/(x^2+y^2)$ on~$\R^2\sm\{0\}$.

In fact the winding number is simply the {\it topological
degree} of $\ga$. Thus it is actually well-defined for
$\ga$ only {\it continuous}, and is invariant under
{\it continuous deformations} of~$\ga$.
\label{un6def1}
\end{dfn}

The motivation for our results is the following theorem
from elementary complex analysis:

\begin{thm} Let\/ $S$ be a domain in $\C$, and 
suppose $f:S\ra\C$ is a holomorphic function, with\/ $f\ne 0$
on $\pd S$. Then the number of zeroes of\/ $f$ in $S^\circ$,
counted with multiplicity, is equal to the winding number
of\/ $f\vert_{\pd S}$ about\/ $0$ along~$\pd S$.
\label{un6thm1}
\end{thm}

As \eq{un4eq5} is a nonlinear version of the Cauchy--Riemann
equations for holomorphic functions, it is natural to expect
that similar results should hold for solutions of \eq{un4eq5}.
We will prove such results.

\subsection{Winding number results for solutions of \eq{un4eq5}}
\label{un61}

Rather than considering with a single solution $u,v$ of \eq{un4eq5},
we shall get more general results by working with two solutions
$u_1,v_1$ and $u_2,v_2$, and treating $(u_1,v_1)-(u_2,v_2)$ like
a holomorphic function for which we wish to count the zeroes. Here is 
the definition of the multiplicity of a zero of~$(u_1,v_1)-(u_2,v_2)$. 

\begin{dfn} Let $S$ be a domain in $\R^2$, let $a\ne 0$, and 
suppose $(u_1,v_1)$ and $(u_2,v_2)$ are solutions of \eq{un4eq5} in
$C^1(S)$. Let $k\ge 1$ be an integer and $(b,c)\in S^\circ$. We say 
that $(u_1,v_1)-(u_2,v_2)$ {\it has a zero of multiplicity $k$ at\/} 
$(b,c)$ if $\pd^ju_1(b,c)=\pd^ju_2(b,c)$ and $\pd^jv_1(b,c)=\pd^j
v_2(b,c)$ for $j=0,\ldots,k-1$, but $\pd^ku_1(b,c)\ne\pd^ku_2(b,c)$ 
and $\pd^kv_1(b,c)\ne\pd^kv_2(b,c)$. Here $\pd$ is the vector 
operator~$(\frac{\pd}{\pd x},\frac{\pd}{\pd y})$.
\label{un6def2}
\end{dfn}

The following lemma justifies this definition, by showing that
every zero of $(u_1,v_1)-(u_2,v_2)$ has a unique multiplicity.

\begin{lem} Let\/ $S$ be a domain in $\R^2$, let\/ $a\ne 0$, and 
let\/ $(u_1,v_1)$ and\/ $(u_2,v_2)$ be solutions of\/ \eq{un4eq5} in
$C^1(S)$, with\/ $(u_1,v_1)\not\equiv(u_2,v_2)$. Suppose $(b,c)\in 
S^\circ$ with\/ $u_1(b,c)=u_2(b,c)$ and\/ $v_1(b,c)=v_2(b,c)$. Then 
$(u_1,v_1)-(u_2,v_2)$ has a zero of multiplicity $k$ at\/ $(b,c)$ for
some unique~$k$.
\label{un6lem1}
\end{lem}

\begin{proof} Since $(u_1,v_1)=(u_2,v_2)$ at one point and
$\pd u_j$ determines $\pd v_j$ by \eq{un4eq5}, it is easy
to see that $u_1\equiv u_2$ if and only if $v_1\equiv v_2$.
But $(u_1,v_1)\not\equiv(u_2,v_2)$ by assumption. Thus
$u_1\not\equiv u_2$ and $v_1\not\equiv v_2$. By Proposition 
\ref{un4prop3}, $u_1,v_1$ and $u_2,v_2$ are real analytic in 
$S^\circ$, and so they are locally given by their Taylor series 
at $(b,c)$. Thus, if $\pd^ju_1(b,c)=\pd^ju_2(b,c)$ for all 
$j=0,1,2,\ldots$ then $u_1\equiv u_2$, a contradiction.

Hence, there exists a unique integer $k\ge 1$ such that 
$\pd^ju_1(b,c)=\pd^ju_2(b,c)$ for $j=0,\ldots,k-1$, and
$\pd^ku_1(b,c)\neq\pd^ku_2(b,c)$. Similarly, there exists 
a unique $l\ge 1$ such that $\pd^jv_1(b,c)=\pd^jv_2(b,c)$ for 
$j=0,\ldots,l-1$, and $\pd^lv_1(b,c)\neq\pd^lv_2(b,c)$. But if 
$\pd^ju_1(b,c)=\pd^ju_2(b,c)$ and $\pd^jv_1(b,c)=\pd^jv_2(b,c)$ 
for $j=0,\ldots,m-1$, one can show from \eq{un4eq5} that
$\pd^mu_1(b,c)=\pd^mu_2(b,c)$ if and only if $\pd^mv_1(b,c)=
\pd^mv_2(b,c)$. This implies that $k=l$, and the lemma follows.
\end{proof}

Next we show that near a zero, $(u_1,v_1)-(u_2,v_2)$ can be
modelled by a genuine holomorphic function, to highest order.

\begin{prop} Let\/ $S$ be a domain in $\R^2$, let\/ $a\ne 0$, and 
let\/ $(u_1,v_1)$ and\/ $(u_2,v_2)$ be solutions of\/ \eq{un4eq5} in
$C^1(S)$. Suppose $(u_1,v_1)-(u_2,v_2)$ has a zero of multiplicity 
$k\ge 1$ at\/ $(b,c)$ in $S^\circ$. Then there exists a nonzero 
complex number $C$ such that
\e
\begin{split}
\la u_1(x,y)+iv_1(x,y)=\la u_2(x,y)&+iv_2(x,y)
+C\bigl(\la(x-b)+i(y-c)\bigr)^k\\
&+O\bigl(\md{x-b}^{k+1}+\md{y-c}^{k+1}\bigr),
\end{split}
\label{un6eq1}
\e
where~$\la=\sqrt{2}\bigl(v_1(b,c)^2+c^2+a^2\bigr)^{1/4}$.
\label{un6prop1}
\end{prop}

\begin{proof} Define polynomials $p(x,y)$, $q(x,y)$ of order $k$ by
\begin{align*}
&p(x,y)=\sum_{j=0}^k\frac{(x-b)^j(y-c)^{k-j}}{j!(k-j)!}\cdot
\frac{\pd^k(u_1-u_2)}{\pd x^j\pd y^{k-j}}(b,c)\\
\text{and}\quad
&q(x,y)=\sum_{j=0}^k\frac{(x-b)^j(y-c)^{k-j}}{j!(k-j)!}\cdot
\frac{\pd^k(v_1-v_2)}{\pd x^j\pd y^{k-j}}(b,c).
\end{align*}
Then as $(u_1,v_1)-(u_2,v_2)$ has a zero of multiplicity $k$ at $(b,c)$,
we see that $p,q$ are nonzero and
\e
\begin{split}
u_1(x,y)&=u_2(x,y)+p(x,y)+O\bigl(\md{x-b}^{k+1}+\md{y-c}^{k+1}\bigr),\\
v_1(x,y)&=v_2(x,y)+q(x,y)+O\bigl(\md{x-b}^{k+1}+\md{y-c}^{k+1}\bigr).
\end{split}
\label{un6eq2}
\e

Taking the difference of equation \eq{un4eq5} for $u_1,v_1$ and 
$u_2,v_2$, the highest order terms at $(b,c)$ imply that
\begin{equation*}
\frac{\pd p}{\pd x}=\frac{\pd q}{\pd y}\quad\text{and}\quad
\frac{\pd q}{\pd x}=-\la^2\frac{\pd p}{\pd y}.
\end{equation*}
But these are the Cauchy--Riemann equations for $\la p+iq$
to be a holomorphic function of $\la x+iy$. Since $p,q$ are 
homogeneous of order $k$ in $(x-b),(y-c)$ it follows that
$\la p(x,y)+iq(x,y)=C\bigl(\la(x-b)+i(y-c)\bigr)^k$ for some 
$C\in\C$, which is nonzero as $p,q$ are nonzero. Combining 
this with \eq{un6eq2} gives~\eq{un6eq1}.
\end{proof}

From \eq{un6eq1} we see that if $(x,y)$ is close to $(b,c)$
in $S^\circ$ but not equal to it then $(u_1,v_1)\ne(u_2,v_2)$
at $(x,y)$. This proves:

\begin{cor} Let\/ $S$ be a domain in $\R^2$, let\/ $a\ne 0$, and 
let\/ $(u_1,v_1)$ and\/ $(u_2,v_2)$ be solutions of\/ \eq{un4eq5} in
$C^1(S)$, with\/ $(u_1,v_1)\not\equiv(u_2,v_2)$. Then the zeroes
of\/ $(u_1,v_1)-(u_2,v_2)$ are isolated in $S^\circ$, that is,
they have no limit points in~$S^\circ$. 

Hence, if\/ $(u_1,v_1)\ne(u_2,v_2)$ at every point of\/ $\pd S$, then
$(u_1,v_1)-(u_2,v_2)$ has finitely many zeroes in~$S$.
\label{un6cor}
\end{cor}

The last part follows because $S$ is compact, and the set of zeroes
of $(u_1,v_1)-(u_2,v_2)$ in $S$ has no limit points. Here is the
main result of this section.

\begin{thm} Let\/ $S$ be a domain in $\R^2$ and\/ $(u_1,v_1)$,
$(u_2,v_2)$ solutions of\/ \eq{un4eq5} in $C^1(S)$ for some
$a\ne 0$, with\/ $(u_1,v_1)\neq(u_2,v_2)$ at every point of\/
$\pd S$. Then $(u_1,v_1)-(u_2,v_2)$ has finitely many zeroes in
$S$. Let there be $n$ zeroes, with multiplicities $k_1,\ldots,k_n$.
Then the winding number of\/ $(u_1,v_1)-(u_2,v_2)$ about\/ $0$
along $\pd S$ is~$\sum_{i=1}^nk_i$.
\label{un6thm2}
\end{thm}

\begin{proof} Let $B_\ep(x,y)$ denote the open ball of radius $\ep$ 
about $(x,y)$ in $\R^2$, and $\,\overline{\!B}_\ep(x,y)$ its 
closure. Let $\ga_\ep(x,y)$ be the circle of radius $\ep$ about 
$(x,y)$, with the natural orientation, and $\bar\ga_\ep(x,y)$ 
the same circle with the reverse orientation.

By Corollary \ref{un6cor} there are finitely many zeroes of 
$(u_1,v_1)-(u_2,v_2)$ in $S$. Let these be $(b_1,c_1),\ldots,(b_n,c_n)$, 
with multiplicities $k_1,\ldots,k_n$ respectively. From \eq{un6eq1} we
see that if $\ep>0$ is sufficiently small then the winding number of
$(u_1,v_1)-(u_2,v_2)$ about 0 along $\ga_\ep(b_i,c_i)$ is $k_i$.
Choose $\ep_1,\ldots,\ep_n>0$ small enough that:
\begin{itemize}
\setlength{\itemsep}{0pt}
\setlength{\parsep}{0pt}
\item $\,\overline{\!B}_{\ep_i}(b_i,c_i)$ lies
in $S^\circ$ for all $i=1,\ldots,n$;
\item $\,\overline{\!B}_{\ep_i}(b_i,c_i)\cap
\,\overline{\!B}_{\ep_k}(b_k,c_k)=\emptyset$ for all $1\le j<k\le n$; and
\item the winding number of $(u_1,v_1)-(u_2,v_2)$ about $0$ along 
$\ga_{\ep_i}(b_i,c_i)$ is~$k_i$.
\end{itemize}

Define $T=S\sm\bigcup_{i=1}^nB_{\ep_i}(b_i,c_i)$. 
Then $(u_1,v_1)-(u_2,v_2)$ has no zeroes in $T$. It follows that the 
winding number of $(u_1,v_1)-(u_2,v_2)$ about 0 along $\pd T$ is zero. 
This can be proved from the definition using Stokes' Theorem, as 
$(u_1-u_2,v_1-v_2)^*(\d\th)$ is a closed 1-form on $T$, 
so~$\int_{\pd T}(u_1-u_2,v_1-v_2)^*(\d\th)=0$. 

Now $\pd T$ is the disjoint union of $\pd S$ and $\bar\ga_{\ep_i}(b_i,c_i)$
for $i=1,\ldots,n$. Thus the winding number of $(u_1,v_1)-(u_2,v_2)$ about
0 along $\pd T$ is the sum of its winding numbers along $\pd S$ and 
$\bar\ga_{\ep_i}(b_i,c_i)$ for $i=1,\ldots,n$. But the winding number 
along $\bar\ga_{\ep_i}(b_i,c_i)$ is $-k_i$, as the winding number 
along $\ga_{\ep_i}(b_i,c_i)$ is $k_i$. Hence the winding number of 
$(u_1,v_1)-(u_2,v_2)$ about 0 along $\pd S$ minus the sum of 
$k_1,\ldots,k_n$ is zero, as we want.
\end{proof}

\subsection{Inverse solutions}
\label{un62}

Recall from \S\ref{un4} that equation \eq{un4eq5} was derived by
beginning with a $\U(1)$-invariant SL 3-fold $N$, and defining
functions $x,y,u$ and $v$ on $N$ by
\begin{equation*}
x=\Re(z_3),\quad u=\Im(z_3), \quad y=\Im(z_1z_2)
\quad\text{and}\quad v=\Re(z_1z_2)
\end{equation*}
for each $(z_1,z_2,z_3)$ in $N$, which also satisfies 
$\ms{z_1}-\ms{z_2}=2a$. Locally we can regard $u,v$ as functions
of $x,y$ (except at branch points), and then the condition that
$N$ be special Lagrangian is equivalent to \eq{un4eq5}.

Consider the map $\si:\C^3\ra\C^3$ defined by $\si(z_1,z_2,z_3)=
(\bar z_1,i\bar z_2,i\bar z_3)$. This is an isometry with
$\si^*(\Re\Om)=-\Re\Om$, and therefore takes SL 3-folds to
SL 3-folds, reversing orientation, as SL 3-folds are calibrated
w.r.t.\ $\Re\Om$. Also, $\si^*(x)=u$, $\si^*(u)=x$, $\si^*(y)=v$ 
and $\si^*(v)=y$, so that $\si$ swaps round $(x,y)$ and $(u,v)$,
and $\si$ preserves the equation~$\ms{z_1}-\ms{z_2}=2a$. 

Therefore, if we regard the SL 3-fold $N$ as a kind of graph
of the function $(x,y)\mapsto(u,v)$, the SL 3-fold $\si(N)$
is the `graph' of the inverse function $(u,v)\mapsto(x,y)$. By 
Proposition \ref{un4prop1}, it follows that $(u,v)$ satisfies
\eq{un4eq5} if and only if its inverse satisfies \eq{un4eq5},
provided a differentiable inverse exists. So we have proved:

\begin{prop} Let\/ $S$ be a domain in $\R^2$, let\/ $a\ne 0$, 
and let\/ $u,v\in C^1(S)$ satisfy \eq{un4eq5}. Define $T=(u,v)[S]$,
and suppose $(u,v):S\ra T$ has a differentiable inverse $(u',v'):
T\ra S$, for $u',v'\in C^1(T)$. Then $u',v'$ satisfy \eq{un4eq5},
and the $\U(1)$-invariant SL\/ $3$-folds $N,N'$ in $\C^3$
corresponding to $u,v$ and\/ $u',v'$ are related by the
involution~$(z_1,z_2,z_3)\mapsto(\bar z_1,i\bar z_2,i\bar z_3)$.
\label{un6prop2}
\end{prop}

One can also easily prove the proposition directly, by expressing
the derivatives of $u',v'$ in terms of those of $u,v$ by matrix
inversion, and observing that \eq{un4eq5} for $u,v$ is equivalent
to \eq{un4eq5} for $u',v'$. We can interpret the proposition as 
an analogue of the fact that the inverses of holomorphic functions
are holomorphic.

\subsection{Nonexistence of $u,v$ with given $u,v,\frac{\pd v}{\pd x},
\frac{\pd v}{\pd y}$ at $(x_0,y_0)$}
\label{un63}

We shall use the `winding number' results of \S\ref{un61} and
the `inverse solution' idea of \S\ref{un62} to show that when
$S,T$ are domains in $\R^2$ and $(\hat u,\hat v):S\ra T$ is a 
solution of \eq{un4eq5}, then maps $(u,v):T\ra S^\circ$ 
satisfying \eq{un4eq5} cannot have certain values of 
$u,v,\frac{\pd v}{\pd x},\frac{\pd v}{\pd y}$ at points 
$(x_0,y_0)$ in~$T^\circ$.

\begin{thm} Let\/ $S,T$ be domains in $\R^2$. Let\/ $a\ne 0$,
$(\hat x_0,\hat y_0)\in S^\circ$, $(\hat u_0,\hat v_0)\in T^\circ$,
and\/ $(\hat p_0,\hat q_0)\in\R^2\sm\{0\}$. Suppose $(\hat u,\hat v):
S\ra T$ is $C^1$ and satisfies \eq{un4eq5} and
\e
\hat u(\hat x_0,\hat y_0)=\hat u_0,\;\> 
\hat v(\hat x_0,\hat y_0)=\hat v_0,\;\>
\frac{\pd\hat v}{\pd x}(\hat x_0,\hat y_0)=\hat p_0 \;\>\text{and\/}\;\>
\frac{\pd\hat v}{\pd y}(\hat x_0,\hat y_0)=\hat q_0.
\label{un6eq3}
\e
Define 
\e
\begin{gathered}
x_0=\hat u_0,\qquad\quad y_0=\hat v_0,\qquad\quad 
u_0=\hat x_0,\qquad\quad v_0=\hat y_0,\\
p_0=-\,\frac{\hat p_0}{\ha(\hat v_0^2\!+\!\hat y_0^2\!+\!a^2)^{-1/2}
\hat p_0^2\!+\!\hat q_0^2}
\;\>\text{and\/}\;\>
q_0=\frac{\hat q_0}{\ha(\hat v_0^2\!+\!\hat y_0^2\!+\!a^2)^{-1/2}
\hat p_0^2\!+\!\hat q_0^2}.
\end{gathered}
\label{un6eq4}
\e
Then there does not exist\/ $(u,v):T\ra S^\circ$ which is $C^1$ 
and satisfies \eq{un4eq5} and
\e
u(x_0,y_0)=u_0, \;\> v(x_0,y_0)=v_0, \;\>
\frac{\pd v}{\pd x}(x_0,y_0)=p_0 \;\>\text{and\/}\;\>
\frac{\pd v}{\pd y}(x_0,y_0)=q_0.
\label{un6eq5}
\e
\label{un6thm3}
\end{thm}

\begin{proof} Suppose for a contradiction that there exists $(u,v):T
\ra S^\circ$ which is $C^1$ and satisfies \eq{un4eq5} and \eq{un6eq5}.
Suppose also that $(\hat u,\hat v):S\ra T$ is injective with nowhere 
vanishing first derivatives. Define $U=(\hat u,\hat v)(S)$. Then $U$ 
is a domain in $\R^2$, and $(\hat u,\hat v):S\ra U$ is an invertible 
map with differentiable inverse.

Let $(u',v'):U\ra S$ be the inverse map. Then by Proposition
\ref{un6prop2}, $u',v'$ satisfy \eq{un4eq5}. As $(\hat u,\hat v)
(\hat x_0,\hat y_0)=(\hat u_0,\hat v_0)$ we see that $(u',v')(x_0,y_0)
=(u_0,v_0)$. Also, since $\hat u,\hat v$ satisfy \eq{un4eq5} we deduce 
from \eq{un6eq3} that
\begin{equation*}
\begin{pmatrix} \frac{\pd\hat u}{\pd x} & \frac{\pd\hat u}{\pd y} \\
\frac{\pd\hat v}{\pd x} & \frac{\pd\hat v}{\pd y} \end{pmatrix}
(\hat x_0,\hat y_0)=\begin{pmatrix} \hat q_0 & 
-\ha(\hat v_0^2+\hat y_0^2+a^2)^{-1/2}\hat p_0 \\
\hat p_0 & \hat q_0\end{pmatrix}.
\end{equation*}
But as $(u',v')$ is the inverse map of $(\hat u,\hat v)$ and
$(\hat u,\hat v)(\hat x_0,\hat y_0)=(x_0,y_0)$ we have
\begin{equation*}
\begin{pmatrix} \frac{\pd u'}{\pd x} & \frac{\pd u'}{\pd y} \\
\frac{\pd v'}{\pd x} & \frac{\pd v'}{\pd y} \end{pmatrix}
(x_0,y_0)=
\begin{pmatrix} \frac{\pd\hat u}{\pd x} & \frac{\pd\hat u}{\pd y} \\
\frac{\pd\hat v}{\pd x} & \frac{\pd\hat v}{\pd y} \end{pmatrix}^{-1}
(\hat x_0,\hat y_0).
\end{equation*}

Combining the last two equations and \eq{un6eq4} shows that
\begin{equation*}
\begin{pmatrix} \frac{\pd u'}{\pd x} & \frac{\pd u'}{\pd y} \\
\frac{\pd v'}{\pd x} & \frac{\pd v'}{\pd y} \end{pmatrix}
(x_0,y_0)=\begin{pmatrix} q_0 & -\ha(v_0^2+y_0^2+a^2)^{-1/2}p_0 \\
p_0 & q_0\end{pmatrix}.
\end{equation*}
Comparing this with \eq{un6eq5} and remembering that $u,v$ satisfy 
\eq{un4eq5}, we see that at $(x_0,y_0)$ we have
\begin{equation*}
\ts u=u',\;\> v=v',\;\>
\frac{\pd u}{\pd x}=\frac{\pd u'}{\pd x},\;\>
\frac{\pd u}{\pd y}=\frac{\pd u'}{\pd y},\;\>
\frac{\pd v}{\pd x}=\frac{\pd v'}{\pd x}\;\>\text{and}\;\>
\frac{\pd v}{\pd y}=\frac{\pd v'}{\pd y}.
\end{equation*}
Thus, $(u',v')-(u,v)$ has a zero of multiplicity at least 2
at $(x_0,y_0)$, in the sense of Definition~\ref{un6def2}.

As $U=(\hat u,\hat v)(S)$ and $(\hat u,\hat v):S\ra T$ we see that 
$U\subseteq T$. Therefore $(u,v)$ and $(u',v')$ are both solutions
of \eq{un4eq5} on the domain $U$. Since $(u',v')$ is an 
orientation-preserving diffeomorphism $U\ra S$, it takes $\pd U$ 
to $\pd S$, and $(u',v')\vert_{\pd U}$ winds once round $\pd S$ in 
the positive sense.

Now $(u,v)$ maps to $S^\circ$ by assumption, and $S$ is contractible.
Therefore the winding number of $(u',v')-(u,v)$ about 0 along $\pd U$
is 1. So by Theorem \ref{un6thm2} the sum of the zeroes of $(u',v')-
(u,v)$ in $U^\circ$, counted with multiplicity, is 1. However, we have
already shown that $(u',v')-(u,v)$ has a zero of multiplicity at least
2 at $(x_0,y_0)$, and $(x_0,y_0)\in U^\circ$ as $(u_0,v_0)\in S^\circ$,
a contradiction.

This proves the theorem under the additional assumption that 
$(\hat u,\hat v):S\ra T$ is injective with nowhere vanishing 
first derivatives. To complete the proof we need to explain how to 
remove this assumption. We can do this using the K\"ahler quotient 
point of view of \S\ref{un42}. Let $\Si$ be the graph of $(u,v)$
in $S\t T$, swapping round the factors $S,T$, and $\hat\Si$ the
graph of $(\hat u,\hat v)$ in~$S\t T$.

We can naturally identify $S\t T$ with a subset of the
K\"ahler quotient $M_a$ discussed in \S\ref{un42}. Thus,
$S\t T$ carries an almost complex structure $J$. Since
$\Si,\hat\Si$ are both quotients of $\U(1)$-invariant SL 3-folds 
in $\C^3$, from \S\ref{un42} we see that $\Si,\hat\Si$ are 
{\it pseudo-holomorphic curves} with respect to~$J$.

Now $\pd\Si\subset S^\circ\t\pd T$ and $\pd\hat\Si\subset\pd S\t T$, and 
$\pd\Si,\pd\hat\Si$ wind once round $\pd T$ and $\pd S$ respectively. 
Therefore the algebraic intersection number $\Si\cap\hat\Si$ is 1. By 
properties of pseudo-holomorphic curves it follows that $\Si,\hat\Si$ 
intersect at only one point, with multiplicity 1. However, the argument 
above shows that $\Si,\hat\Si$ intersect with multiplicity at least 2 
at $(u_0,v_0,x_0,y_0)$, a contradiction, and the theorem is complete.
\end{proof}

This theorem will be used in \cite{Joyc6} to construct a priori
estimates for $\frac{\pd u}{\pd x},\frac{\pd u}{\pd y},\frac{\pd
v}{\pd x}$ and $\frac{\pd v}{\pd y}$ for bounded solutions $u,v$
of~\eq{un4eq5}.

\section{Rewriting \eq{un4eq5} in terms of a potential $f$}
\label{un7}

Let $S$ be a domain in $\R^2$, as in Definition \ref{un3def1}, 
and fix $a\ne 0$ in $\R$. We shall study differentiable 
functions $u,v:S\ra\R$ satisfying equation \eq{un4eq5} in $S$, 
and also certain {\it boundary conditions} on $\pd S$. As 
$\frac{\pd u}{\pd x}=\frac{\pd v}{\pd y}$, we can write $u,v$ in 
terms of a {\it potential}\/ $f:S\ra\R$ with $u=\frac{\pd f}{\pd y}$
and~$v=\frac{\pd f}{\pd x}$. 

\begin{prop} Let\/ $S$ be a domain in $\R^2$ and\/ $u,v\in C^1(S)$
satisfy \eq{un4eq5} for $a\ne 0$. Then there exists $f\in C^2(S)$
with\/ $\frac{\pd f}{\pd y}=u$, $\frac{\pd f}{\pd x}=v$ and
\e
P(f)=\Bigl(\Bigl(\frac{\pd f}{\pd x}\Bigr)^2+y^2+a^2
\Bigr)^{-1/2}\frac{\pd^2f}{\pd x^2}+2\,\frac{\pd^2f}{\pd y^2}=0.
\label{un7eq1}
\e
This $f$ is unique up to addition of a constant, $f\mapsto f+c$.
Conversely, all solutions of\/ \eq{un7eq1} yield solutions 
of\/~\eq{un4eq5}. 
\label{un7prop1}
\end{prop}

\begin{proof} Define a 1-form $\al$ on $S$ by $\al=v(x,y)\d x+u(x,y)\d y$.
Then $\d\al=0$ as $\frac{\pd v}{\pd y}=\frac{\pd u}{\pd x}$, so $\al$ is
closed. As $S$ is contractible, $\al$ is exact, and so $\al=\d f$ for
some $f\in C^2(S)$, unique up to addition of a constant. Equating
coefficients of $\d x$ and $\d y$ in $\al=\d f$ gives
$\frac{\pd f}{\pd x}=v$, $\frac{\pd f}{\pd y}=u$. Equation \eq{un7eq1}
follows by substituting these into the first equation of \eq{un4eq5}
and multiplying by $(v^2+y^2+a^2)^{-1/2}$. The converse is easy.
\end{proof}

Now \eq{un7eq1} is a {\it second-order quasilinear elliptic equation},
so Theorem \ref{un3thm1}~gives:

\begin{thm} Let\/ $S$ be a domain in $\R^2$, let\/ $a\ne 0$,
and suppose $f\in C^2(S)$ satisfies \eq{un7eq1} with\/ 
$f\vert_{\pd S}=\phi\in C^2(\pd S)$. Then $f$ is real analytic
in $S^\circ$, and if\/ $\phi\in C^{k+2,\al}(\pd S)$ for $k\ge 0$
and\/ $\al\in(0,1)$ then $f\in C^{k+2,\al}(S)$, and if\/
$\phi\in C^\iy(\pd S)$ then~$f\in C^\iy(S)$.
\label{un7thm1}
\end{thm}

As \eq{un7eq1} is of the form \eq{un3eq1} with $b^i\equiv c\equiv 0$, 
by the maximum principle, Theorem \ref{un3thm2}, we deduce:

\begin{lem} Let\/ $S$ be a domain in $\R^2$, let\/ $a\ne 0$, and
suppose $f\in C^2(S)$ is a solution of\/ \eq{un7eq1}. Then the
maximum and minimum of\/ $f$ are achieved on~$\pd S$.
\label{un7lem1}
\end{lem}

Equation \eq{un7eq1} may also be written
\e
P(f)=\frac{\pd}{\pd x}\Bigl[A\Bigl(y,\frac{\pd f}{\pd x}\Bigr)\Bigr]
+2\,\frac{\pd^2f}{\pd y^2}=0,
\label{un7eq2}
\e
where $A(y,v)$ is defined to be
\e
A(y,v)=\int_0^v\bigl(w^2+y^2+a^2\bigr)^{-1/2}\,\d w,
\;\>\text{so that}\;\>
\frac{\pd A}{\pd v}=\bigl(v^2+y^2+a^2\bigr)^{-1/2}.
\label{un7eq3}
\e
Equation \eq{un7eq2} is equivalent to \eq{un7eq1}, but is in
{\it divergence form}.

Calculation shows that we may write $A$ explicitly as
\begin{equation*}
A(y,v)=\log\biggl[\frac{(v^2+y^2+a^2)^{1/2}+v}{(y^2+a^2)^{1/2}}\biggr]
=\log\biggl[\frac{(y^2+a^2)^{1/2}}{(v^2+y^2+a^2)^{1/2}-v}\biggr].
\end{equation*}
Note that $A$ is undefined when $a=y=0$. That is, if $a=0$ then
$A$ is undefined along the $x$-axis. 

\subsection{Expressing \eq{un7eq1} as an Euler--Lagrange equation}
\label{un71}

We shall show that equation \eq{un7eq1} is in fact the
{\it Euler--Lagrange equation} of a certain functional 
$I:C^{0,1}(S)\ra\R$. Fix $a\ne 0$, and define a function 
$B(y,v)$ by $B(y,v)=\int_0^vA(w,y)\d w$, so that 
$\frac{\pd B}{\pd v}(y,v)=A(y,v)$. Define a function $F$
on $S\t\R^2$ by
\e
F(x,y,u,v)=B(y,v)+u^2,
\label{un7eq4}
\e
and define a functional $I:C^{0,1}(S)\ra\R$ by
\e
I(f)=\int_SF\Bigl(x,y,\frac{\pd f}{\pd x},\frac{\pd f}{\pd y}\Bigr)\d x\,\d y.
\label{un7eq5}
\e
The {\it Euler--Lagrange} equation for $I$ is
\begin{equation*}
\frac{\pd}{\pd x}\biggl[\frac{\pd F}{\pd v}\Bigl(x,y,
\frac{\pd f}{\pd x},\frac{\pd f}{\pd y}\Bigr)\biggr]+
\frac{\pd}{\pd y}\biggl[\frac{\pd F}{\pd u}\Bigl(x,y,
\frac{\pd f}{\pd x},\frac{\pd f}{\pd y}\Bigr)\biggr]=0.
\end{equation*}
From \eq{un7eq4} this becomes
\begin{equation*}
\frac{\pd}{\pd x}\biggl[\frac{\pd B}{\pd v}\Bigl(y,\frac{\pd f}{\pd x}
\Bigr)\biggr]+
\frac{\pd}{\pd y}\biggl[2\,\frac{\pd f}{\pd y}\biggr]=0,
\end{equation*}
and this is equivalent to \eq{un7eq2}, since $\frac{\pd B}{\pd v}(y,v)
=A(y,v)$. Thus we have proved:

\begin{prop} Equations \eq{un7eq1} and\/ \eq{un7eq2} are equivalent to 
the Euler--Lagrange equation of the functional\/ $I:C^{0,1}(S)\ra\R$ 
defined in~\eq{un7eq5}.
\label{un7prop2}
\end{prop}

We could use this to solve the Dirichlet problem for \eq{un7eq1}
on $S$, by choosing a minimizing sequence $(f_n)_{n=1}^\iy$ for
$I$ amongst all $f\in C^{0,1}(S)$ with $f\vert_{\pd S}=\phi$ for
some $\phi\in C^{k+2,\al}(\pd S)$, and then showing that $f_n$
converges to a solution as $n\ra\iy$. But we will instead do it
by more elementary methods in~\S\ref{un73}.

\subsection{Super- and subsolutions of \eq{un7eq1}}
\label{un72}

Let $S$ be a domain in $\R^2$, let $a\ne 0$, and let $P$ be
the operator defined in \eq{un7eq1}. A function $f\in C^2(S)$ 
is called a {\it supersolution} of \eq{un7eq1} if $P(f)\ge 0$
on $S$, and a {\it subsolution} if $P(f)\le 0$. Sub- and
supersolutions $f,f'$ with $f\le f'$ on $\pd S$ satisfy
$f\le f'$ on~$S$.

\begin{prop} Let\/ $S$ be a domain in $\R^2$ and\/ $a\ne 0$.
Suppose $f,f'\in C^2(S)$ satisfy $P(f)\ge 0$ on $S$ and\/
$P(f')\le 0$ on $S$, where $P$ is defined in \eq{un7eq1}.
If\/ $f\le f'$ on $\pd S$ then $f\le f'$ on $S$, and if\/
$f<f'$ on $\pd S$ then $f<f'$ on~$S$.
\label{un7prop3}
\end{prop}

\begin{proof} Applying the Mean Value Theorem to $F(z)=(z^2+y^2+a^2
)^{-1/2}$ on the interval $[\frac{\pd f}{\pd x},\frac{\pd f'}{\pd x}]$
we find that
\begin{equation*}
\ts\bigl((\frac{\pd f}{\pd x})^2\!+\!y^2\!+\!a^2\bigr)^{-1/2}\!-\!
\bigl((\frac{\pd f'}{\pd x})^2\!+\!y^2\!+\!a^2\bigr)^{-1/2}\!=\!
-w\bigl(w^2\!+\!y^2\!+\!a^2\bigr)^{-3/2}
\bigl(\frac{\pd f}{\pd x}\!-\!\frac{\pd f'}{\pd x}\bigr)
\end{equation*}
for some $w$ between $\frac{\pd f}{\pd x}$ and $\frac{\pd f'}{\pd x}$.
Using \eq{un7eq1} to expand $P(f)-P(f')$ and rearranging then gives
\begin{align*}
P(f)-P(f')=&\bigl(({\ts\frac{\pd f}{\pd x}})^2+y^2+a^2\bigr)^{-1/2}
\frac{\pd^2}{\pd x^2}\bigl(f-f'\bigr)
+2\frac{\pd^2}{\pd y^2}\bigl(f-f'\bigr)\\
&-\bigl(w(w^2+y^2+a^2)^{-3/2}{\ts\frac{\pd^2f'}{\pd x^2}}\bigr)
\frac{\pd}{\pd x}\bigl(f-f'\bigr).
\end{align*}

We may regard the right hand side of this equation as $L(f-f')$,
where $L$ is a {\it linear elliptic operator} of the form \eq{un3eq1}
with $c\equiv 0$. Thus $L(f-f')=P(f)-P(f')\ge 0$, as $P(f)\ge 0$
and $P(f')\le 0$. So by the maximum principle, Theorem \ref{un3thm2},
the maximum of $f-f'$ is achieved on $\pd S$. Thus, if $f-f'\le 0$
on $\pd S$ then $f-f'\le 0$ on $S$, and if $f-f'<0$ on $\pd S$ then
$f-f'<0$ on~$S$.
\end{proof}

In particular, if $P(f)=P(f')=0$ and $f\vert_{\pd S}=f'\vert_{\pd S}$,
then the proposition implies that $f\le f'$ and (exchanging $f,f'$)
that $f'\le f$, so that $f=f'$. This implies uniqueness of solutions
of the Dirichlet problem for~\eq{un7eq1}.

\subsection{The Dirichlet problem for $f$}
\label{un73}

Observe that \eq{un7eq1} is of the form \eq{un3eq5}. Therefore 
Theorem \ref{un3thm3} applies to give existence for the
Dirichlet problem for $f$, and an a priori bound for $\cnm{f}{1}$.
Combining this with the real analyticity in Theorem \ref{un7thm1}
and the uniqueness following from Proposition \ref{un7prop3} gives:

\begin{thm} Let\/ $S$ be a strictly convex domain in $\R^2$, 
and let\/ $a\ne 0$, $k\ge 0$ and\/ $\al\in(0,1)$. Then for 
each\/ $\phi\in C^{k+2,\al}(\pd S)$ there exists a unique 
solution $f$ of\/ \eq{un7eq1} in $C^{k+2,\al}(S)$ with\/ 
$f\vert_{\pd S}=\phi$. This $f$ is real analytic in $S^\circ$, 
and satisfies $\cnm{f}{1}\le C\cnm{\phi}{2}$, for some $C>0$
depending only on~$S$.
\label{un7thm2}
\end{thm}

Thus, the Dirichlet problem for \eq{un7eq1} is uniquely solvable
in a strictly convex domain. Combining the theorem with Propositions 
\ref{un4prop1} and \ref{un7prop1}, we get an existence and uniqueness 
result for $\U(1)$-invariant special Lagrangian 3-folds in $\C^3$ 
satisfying certain boundary conditions.

However, solving the Dirichlet problem in a general, nonconvex domain 
is more difficult, as to get an a priori estimate for $\md{\pd f}$ on 
$\pd S$ one needs to find super- and subsolutions of \eq{un7eq1} 
satisfying certain equalities and inequalities on $\pd S$, and this 
does not seem easy to do in an elementary way. The point about 
strictly convex domains is that one can use affine functions as 
super- and subsolutions to estimate~$\md{\pd f}$.

An analogue of Theorem \ref{un7thm2} for the case $a=0$ will be given
in \cite[Th.~7.1]{Joyc6}, which shows that \eq{un7eq2} has a unique
solution $f\in C^1(S)$ with weak second derivatives, and $f\vert_{\pd
S}=\phi$. But $f$ may have singular points, at which it is only once
differentiable.

By looking closely at the proofs of existence and uniqueness
of $f$, one can show that small changes in $\phi$ and $a$ result
in small changes in $f$, where `small' may be interpreted in the
$C^{k+2,\al}$ sense. Hence we may prove:

\begin{thm} Let\/ $S$ be a strictly convex domain in $\R^2$, $k\ge 0$ 
and\/ $\al\in(0,1)$. Then the map $C^{k+2,\al}(\pd S)\t\bigl(\R\sm\{0\}
\bigr)\ra C^{k+2,\al}(S)$ taking $(\phi,a)\mapsto f$ is continuous, 
where $f$ is the unique solution of\/ \eq{un7eq1} with\/ 
$f\vert_{\pd S}=\phi$ constructed in Theorem~\ref{un7thm2}.
\label{un7thm3}
\end{thm}

Presumably this map $(\phi,a)\mapsto f$ is also smooth. An
extension of Theorem \ref{un7thm3} to include the case $a=0$ is
given in \cite[Th.~7.2]{Joyc6}, but with the $C^1$ rather than
the $C^{k+2,\al}$ topology on~$f$.

\subsection{Winding number results for potentials}
\label{un74}

We shall now extend some of the `winding number' results of 
\S\ref{un6} to the situation of this section. We begin with
a definition.

\begin{dfn} Let $S$ be a domain in $\R^2$, and $\phi\in C^2(\pd S)$.
Choose a smooth parametrization
\e
\R/2\pi\Z\ra\pd S, \quad\text{written}\quad
\th\mapsto\bigl(x(\th),y(\th)\bigr) \quad\text{for $\th\in\R/2\pi\Z$,}
\label{un7eq6}
\e
and regard $\phi$ as a function of $\th$. We call $\phi$ a {\it Morse
function} if $\frac{\d\phi}{\d\th}$ is zero at only finitely many
points in $\pd S$, and $\frac{\d^2\phi}{\d\th^2}$ is nonzero at each
of these points.

It can be shown that this definition is independent of the 
parametrization \eq{un7eq6}, and that the Morse functions are an 
open dense subset of $C^2(\pd S)$. Also, each stationary point of 
$\phi$ on $\pd S$ is either a local maximum or a local minimum, as 
$\frac{\d^2\phi}{\d\th^2}\ne 0$, and there are the same number of 
each, so $\phi$ has exactly $l$ local maxima and $l$ local minima 
for some~$l\ge 1$.
\label{un7def1}
\end{dfn}

If $f\in C^2(S)$ and $f\vert_{\pd S}$ is a Morse function, we
can relate the winding number of $\pd f$ round $\pd S$ to the 
number of local maxima and minima of $f$ on~$\pd S$.

\begin{prop} Let\/ $S$ be a domain in $\R^2$, and\/ $f\in C^2(S)$
with\/ $f\vert_{\pd S}=\phi$. Suppose that $\pd f\ne 0$ at each
point of\/ $\pd S$ and the winding number of\/ $\pd f$ about\/ $0$
along $\pd S$ is $k$, and that\/ $\phi\in C^2(\pd S)$ is a Morse
function with\/ $l$ local maxima and\/ $l$ local minima for some
$l\ge 1$. Then~$1-l\le k\le 1+l$.
\label{un7prop4}
\end{prop}

\begin{proof} Choose a smooth, positively oriented parametrization for 
$\pd S$ as in \eq{un7eq6}. Let the $l$ local maxima of $\phi$ 
be at $\th=\al_j$ and the $l$ local minima at $\th=\be_j$, where
$\al_1,\ldots,\al_l$ and $\be_1,\ldots,\be_l$ lie in $\R/2\pi\Z$,
and are arranged in the cyclic order $\al_1,\be_1,\al_2,\be_2,\ldots,
\be_l,\al_{l+1}=\al_1$. Define $(\ga,\de):\R/2\pi\Z\ra\R^2\sm\{0\}$
by $(\ga,\de)=\frac{\d}{\d\th}\bigl(x(\th),y(\th)\bigr)$, so that
$(\ga,\de)(\th)$ is tangent to $\pd S$ at~$\bigl(x(\th),y(\th)\bigr)$.

Then $\frac{\d\phi}{\d\th}=\ga\frac{\pd f}{\pd x}+\de\frac{\pd f}{\pd y}$
on $\pd S$. Therefore we have
\e
\ga\frac{\pd f}{\pd x}+\de\frac{\pd f}{\pd y}
\begin{cases}=0, & \quad \th=\al_j\;\>\text{or}\;\> \th=\be_j,
\;\> j=1,\ldots,l, \\
<0, & \quad \al_j<\th<\be_j,\quad j=1,\ldots,l,\\
>0, & \quad \be_j<\th<\al_{j+1},\quad j=1,\ldots,l,
\end{cases}
\label{un7eq7}
\e
using the cyclic order on $\R/2\pi\Z$. Also, as $\pd f$ is nonzero
at each point of $\pd S$, we see that if $\th=\al_j$ or $\be_j$
then $\de\frac{\pd f}{\pd x}-\ga\frac{\pd f}{\pd y}\ne 0$. Define
\ea
\eta_j&=\begin{cases} \phantom{-}1, & 
\de(\al_j)\frac{\pd f}{\pd x}\bigl(x(\al_j),y(\al_j)\bigr)-
\ga(\al_j)\frac{\pd f}{\pd y}\bigl(x(\al_j),y(\al_j)\bigr)<0,\\
\phantom{-}0, & 
\de(\al_j)\frac{\pd f}{\pd x}\bigl(x(\al_j),y(\al_j)\bigr)-
\ga(\al_j)\frac{\pd f}{\pd y}\bigl(x(\al_j),y(\al_j)\bigr)>0,
\end{cases}
\label{un7eq8}\\
\text{and}\quad \ze_j&=\begin{cases} -1, & 
\de(\be_j)\frac{\pd f}{\pd x}\bigl(x(\be_j),y(\be_j)\bigr)-
\ga(\be_j)\frac{\pd f}{\pd y}\bigl(x(\be_j),y(\be_j)\bigr)<0,\\
\phantom{-}0, & 
\de(\be_j)\frac{\pd f}{\pd x}\bigl(x(\be_j),y(\be_j)\bigr)-
\ga(\be_j)\frac{\pd f}{\pd y}\bigl(x(\be_j),y(\be_j)\bigr)>0.
\end{cases}
\label{un7eq9}
\ea

Now we can use equations \eq{un7eq7}--\eq{un7eq9} to compare the 
winding numbers of $(\de,-\ga)$ and $\pd f$ about 0 along $\pd S$,
as they tell us when the direction of $\pd f$ crosses that of
$\pm(\de,-\ga)$. But the winding number of $(\de,-\ga)$ about 0 
along $\pd S$ is 1, as it is an outward normal vector to $\pd S$. 
Using this it is easy to show that the winding number of $\pd f$ 
about 0 along $\pd S$ is $k=1+\sum_{j=1}^l\eta_j+\sum_{j=1}^l\ze_j$. 
As $\eta_j$ is 0 or 1 and $\ze_j$ is 0 or $-1$, we see 
that~$1-l\le k\le 1+l$.
\end{proof}

Here is the main result of this subsection:

\begin{thm} Let\/ $S$ be a domain in $\R^2$ and\/ $f_1,f_2\in C^2(S)$
satisfy \eq{un7eq1} for $a\ne 0$ with\/ $f_j\vert_{\pd S}=\phi_j$.
Set\/ $u_j=\frac{\pd f_j}{\pd y}$ and\/ $v_j=\frac{\pd f_j}{\pd x}$, so
that\/ $u_j,v_j\in C^1(S)$ satisfy \eq{un4eq5}. Suppose $\phi_1-\phi_2$
is a Morse function on $\pd S$, with\/ $l$ local maxima and\/ $l$ local
minima. Then $(u_1,v_1)-(u_2,v_2)$ has $n$ zeroes in $S^\circ$ with
multiplicities $k_1,\ldots,k_n$ and\/ $m$ zeroes on $\pd S$,
where~$\sum_{i=1}^nk_i+m\le l-1$.
\label{un7thm4}
\end{thm}

\begin{proof} First suppose, for simplicity, that $(u_1,v_1)
\neq(u_2,v_2)$ at every point of $\pd S$. Then $m=0$, and the
theorem in this case follows from Theorem \ref{un6thm2} and 
Propositions \ref{un7prop1} and \ref{un7prop4}, noting that
\begin{equation*}
\ts\pd(f_1-f_2)=\bigl(\frac{\pd}{\pd x}(f_1-f_2),
\frac{\pd}{\pd y}(f_1-f_2)\bigr)=(v_1-v_2,u_1-u_2),
\end{equation*}
so that the winding number of $\pd(f_1-f_2)$ about 0 along $\pd S$ is 
$-\sum_{i=1}^nk_i$, by Theorem \ref{un6thm2}. It remains to prove the 
result in the case when $(u_1,v_1)=(u_2,v_2)$ at $m\ge 1$ points 
$(x_0,y_0)$ in~$\pd S$. 

Then $\pd(f_1-f_2)=0$ at $(x_0,y_0)$, so $(x_0,y_0)$ must be one of 
the $l$ local maxima or $l$ local minima of $\phi_1-\phi_2$. Thus $m$
is finite. Furthermore, as $\phi_1-\phi_2$ is a Morse function
$\frac{\d^2}{\d\th^2}(\phi_1-\phi_2)\ne 0$ at $(x_0,y_0)$, which implies
that $\pd(u_1,v_1)\ne\pd(u_2,v_2)$ at $(x_0,y_0)$, and therefore
$(x_0,y_0)$ is an {\it isolated}\/ zero of $(u_1,v_1)-(u_2,v_2)$. 
By Corollary \ref{un6cor} and compactness of $S$ we deduce that 
$(u_1,v_1)-(u_2,v_2)$ has finitely many zeroes in $S^\circ$, so we
can suppose there are $n$ zeroes, with multiplicities~$k_1,\ldots,k_n$.

For $\ep\ge 0$, define $S_\ep$ to be the subset of $(x,y)\in S$
with distance at least $\ep$ from $\pd S$, so that $S_0=S$. Choose 
$\ep>0$ sufficiently small that $S_\ep$ is a domain, and $S_\ep^\circ$ 
contains all the $n$ zeroes of $(u_1,v_1)-(u_2,v_2)$ in $S^\circ$, and 
$f\vert_{\pd S_\ep}$ is also a Morse function with $l$ local maxima 
and $l$ local minima. It is easy to see that this is possible. Then
$(u_1,v_1)\neq(u_2,v_2)$ at every point of $\pd S_\ep$, as the zeroes
of $(u_1,v_1)-(u_2,v_2)$ in $S^\circ$ lie in~$S_\ep^\circ$.

Let $k$ be the winding number of $\pd(f_1-f_2)$ about 0 along $\pd 
S_\ep$. Then Proposition \ref{un7prop4} shows that $1-l\le k\le 1+l$. 
However, we can improve the result in this case. Recall that 
$\pd(f_1-f_2)=0$ at $m$ out of the $2l$ local maxima and minima of 
$\phi_1-\phi_2$ on $\pd S$. Using \eq{un7eq1} we can show that if
$\th=\al_j$ is one of these $m$ points then $\eta_j=1$ in \eq{un7eq8}
at the corresponding local maximum in $\pd S_\ep$, and if $\th=\be_j$
is one of the $m$ points then $\ze_j=0$ in \eq{un7eq9} at the 
corresponding local minimum in $\pd S_\ep$. Thus, the proof of 
Proposition \ref{un7prop4} shows that $1-l+m\le k\le 1+l$. But 
applying Theorem \ref{un6thm2} gives $k=-\sum_{i=1}^nk_i$, and 
the theorem follows. 
\end{proof}

The theorem can be used in conjunction with Theorem \ref{un7thm2},
the solution of the Dirichlet problem for $f$ on a strictly convex
domain. In this case, we would know $\phi_1,\phi_2$ explicitly, but
would otherwise know little about the $f_j,u_j$ or $v_j$. The
theorem tells us something about $(u_1,v_1)$ and $(u_2,v_2)$,
using only the boundary data~$\phi_1,\phi_2$.

Using Theorems \ref{un7thm2} and \ref{un7thm3} we can drop the
condition that $\phi_1-\phi_2$ is Morse, requiring instead that
it has only finitely many local maxima and minima.

\begin{thm} Let\/ $S$ be a strictly convex domain in $\R^2$,
let\/ $a\ne 0$, $\al\in(0,1)$, and\/ $f_1,f_2\in C^{2,\al}(S)$
satisfy \eq{un7eq1} with\/ $f_j\vert_{\pd S}=\phi_j$. Set\/
$u_j=\frac{\pd f_j}{\pd y}$ and\/ $v_j=\frac{\pd f_j}{\pd x}$,
so that\/ $u_j,v_j\in C^{1,\al}(S)$ satisfy \eq{un4eq5}.
Suppose $\phi_1-\phi_2$ has exactly $l$ local maxima and\/
$l$ local minima on $\pd S$. Then $(u_1,v_1)-(u_2,v_2)$ has
$n$ zeroes in $S^\circ$ with multiplicities $k_1,\ldots,k_n$,
where~$\sum_{i=1}^nk_i\le l-1$.
\label{un7thm5}
\end{thm}

\begin{proof} It is not difficult to construct a smooth family
$\phi_1^t\in C^{2,\al}(\pd S)$ for $t\in(0,1]$, such that
$\phi_1^t\ra\phi_1$ as $t\ra 0_+$, and $\phi_1^t$ is a Morse
function with $l$ local maxima and $l$ local minima, at the same
points as $\phi_1$. Let $f_1^t$ be the solution of \eq{un7eq1}
given by Theorem \ref{un7thm2} with $f_1^t\vert_{\pd S}=\phi_1^t$,
and set $u_1^t=\frac{\pd f_1^t}{\pd y}$ and~$v_1^t=\frac{\pd
f_1^t}{\pd x}$.

Then the sum of the zeroes of $(u_1^t,v_1^t)-(u_2,v_2)$ in
$S^\circ$ with multiplicity is no more than $l-1$ for all
$t\in(0,1]$, by Theorem \ref{un7thm4}. Also $(u_1^t,v_1^t)
\ra(u_1,v_1)$ in $C^{1,\al}(S)$ as $t\ra 0_+$ by Theorem
\ref{un7thm3}. Combining these using a limiting argument
we find that the sum of the zeroes of $(u_1,v_1)-(u_2,v_2)$
in $S^\circ$ with multiplicity is no more than~$l-1$.
\end{proof}

Note that the theorem does not limit the number of zeroes of
$(u_1,v_1)-(u_2,v_2)$ on $\pd S$, which can appear at any
stationary point of~$\phi_1-\phi_2$.

\section{Another approach to solving \eq{un4eq5}}
\label{un8}

In Proposition \ref{un4prop3} we showed that if $S$ is a domain
in $\R^2$ and $u,v\in C^1(S)$ satisfy \eq{un4eq5} then $v$ satisfies 
\eq{un4eq12} in $S^\circ$. Conversely, if $v\in C^2(S)$ satisfies 
\eq{un4eq12} then using \eq{un4eq5} to find $\frac{\pd u}{\pd x}$,
$\frac{\pd u}{\pd y}$, it is easy to show that as $S$ is contractible
there exists $u\in C^2(S)$, unique up to addition of a constant, such
that $u,v$ satisfy \eq{un4eq5}. In this way we prove:

\begin{prop} Let\/ $S$ be a domain in $\R^2$ and\/ $u,v\in C^2(S)$
satisfy \eq{un4eq5} for $a\ne 0$. Then
\e
Q(v)=\frac{\pd}{\pd x}\Bigl[\bigl(v^2+y^2+a^2\bigr)^{-1/2}
\frac{\pd v}{\pd x}\Bigr]+2\,\frac{\pd^2v}{\pd y^2}=0.
\label{un8eq1}
\e
Conversely, if\/ $v\in C^2(S)$ satisfies \eq{un8eq1} then 
there exists $u\in C^2(S)$, unique up to addition of a 
constant\/ $u\mapsto u+c$, such that $u,v$ satisfy~\eq{un4eq5}.
\label{un8prop1}
\end{prop}

Equation \eq{un8eq1} is a {\it second-order quasilinear elliptic 
equation} upon $v$. It is also in {\it divergence form}. By elliptic 
regularity, Theorem \ref{un3thm1}, we get:

\begin{prop} Let\/ $S$ be a domain in $\R^2$, let\/ $a\ne 0$, 
and suppose $v\in C^2(S)$ is a solution of\/ \eq{un8eq1} with\/ 
$v\vert_{\pd S}=\phi$ for some $\phi\in C^2(\pd S)$. Then $v$ is 
real analytic in $S^\circ$, and if\/ $\phi\in C^{k+2,\al}(\pd S)$
for $k\ge 0$ and\/ $\al\in(0,1)$ then $v\in C^{k+2,\al}(S)$,
and if\/ $\phi\in C^\iy(\pd S)$ then~$v\in C^\iy(S)$.
\label{un8prop2}
\end{prop}

Taking the derivative in \eq{un8eq1} gives the equivalent
\e
\bigl(v^2+y^2+a^2\bigr)^{-1/2}\frac{\pd^2 v}{\pd x^2}
-\frac{v}{\bigl(v^2+y^2+a^2\bigr)^{3/2}}
\Bigl(\frac{\pd v}{\pd x}\Bigr)^2+2\,\frac{\pd^2v}{\pd y^2}=0.
\label{un8eq2}
\e
This is of the form \eq{un3eq1} with $c=0$. Therefore by the 
maximum principle, Theorem \ref{un3thm2}, we have:

\begin{lem} Let\/ $S$ be a domain in $\R^2$, let\/ $a\ne 0$, and
suppose $v\in C^2(S)$ is a solution of\/ \eq{un8eq1}. Then the
maximum and minimum of\/ $v$ are achieved on~$\pd S$.
\label{un8lem1}
\end{lem}

\subsection{Super- and subsolutions of \eq{un8eq1}}
\label{un81}

We now carry out the programme of \S\ref{un72} for
equation~\eq{un8eq1}.

\begin{prop} Let\/ $T$ be a closed, bounded subset of\/ $\R^2$ 
whose boundary $\pd T=T\sm T^\circ$ is a piecewise-smooth closed 
curve, and let\/ $a\ne 0$. Suppose $v,v'\in C^2(T)$ satisfy 
$Q(v)\ge 0$, $Q(v')\le 0$ and\/ $v\ge v'$ on $T$, where $Q$ 
is defined in \eq{un8eq1}, and\/ $v=v'$ on $\pd T$. Then 
$v=v'$ on~$T$.
\label{un8prop3}
\end{prop}

\begin{proof} Choose $C>0$ such that $y^2\le C$ on $T$. Then we have
\e
\begin{split}
0&\ge -\int_T(C-y^2)\bigl[Q(v)-Q(v')\bigr]\d x\,\d y\\
&=-\int_T(C-y^2)\Bigl[\frac{\pd}{\pd x}\Bigl(\bigl(v^2\!+\!y^2\!+\!a^2
\bigr)^{-1/2}\frac{\pd v}{\pd x}-\bigl((v')^2\!+\!y^2\!+\!a^2\bigr)^{-1/2}
\frac{\pd v'}{\pd x}\Bigr)\\
&\qquad\qquad\qquad\qquad
+2\frac{\pd^2}{\pd y^2}\bigl(v-v'\bigr)\Bigr]\d x\,\d y\\
&=\int_{\pd T}(C-y^2)\bigl(v^2\!+\!y^2\!+\!a^2\bigr)^{-1/2}
\Bigl(-\frac{\pd}{\pd x}\bigl(v-v'\bigr)\Bigr)\d y\\
&\qquad
+2\int_{\pd T}(C-y^2)\frac{\pd}{\pd y}\bigl(v-v'\bigr)\d x
+4\int_T(v-v')\d x\,\d y,
\end{split}
\label{un8eq3}
\e
using integration by parts, and the fact that $v=v'$ on~$\pd T$.

We claim that all three integrals on the final line of \eq{un8eq3}
are nonnegative. For the first integral, as $v-v'=0$ on $\pd T$
and $v-v'\ge 0$ on $T$, we see that if $(x,y)\in\pd T$ and $w$
is a vector in $\R^2$ pointing outwards from $T$ at $(x,y)$ then
$\pd_{w}(v-v')\vert_{(x,y)}\le 0$, and if $w$ points inwards
from $T$ then $\pd_{w}(v-v')\vert_{(x,y)}\ge 0$. But $w=
\frac{\pd}{\pd x}$ points outwards from $T$ at $(x,y)$ if and
only if $\d y\vert_{\pd T}$ is a positive 1-form on $\pd T$ at
$(x,y)$, with the natural orientation on~$\pd T$.

Hence $-\frac{\pd}{\pd x}\bigl(v-v'\bigr)\d y\vert_{\pd T}$
is a nonnegative 1-form on $\pd T$, and the first integral on the
final line of \eq{un8eq3} is nonnegative. Similarly $\frac{\pd}{\pd
y}\bigl(v-v'\bigr)\d x\vert_{\pd T}$ is a nonnegative 1-form, so
the second integral is nonnegative, and the third integral is
nonnegative as $v-v'\ge 0$. But the sum of the three is nonpositive
by \eq{un8eq3}. Thus all three integrals are zero, and
$\int_T(v-v')\d x\,\d y=0$. As $v-v'\ge 0$ and $v,v'$ are
continuous, this implies that $v=v'$ on~$T$.
\end{proof}

Using this we can prove an analogue of Proposition \ref{un7prop3}
for \eq{un8eq1}. The restriction to real analytic $v,v'$ is not
really necessary, but simplifies the proof.

\begin{prop} Let\/ $S$ be a domain in $\R^2$ and\/ $a\ne 0$.
Suppose $v,v'\in C^2(S)$ are real analytic in $S^\circ$ and
satisfy $v\le v'$ on $\pd S$, $Q(v)\ge 0$ on $S$ and\/
$Q(v')\le 0$ on $S$, where $Q$ is defined in \eq{un8eq1}.
Then $v\le v'$ on $S$.
\label{un8prop4}
\end{prop}

\begin{proof} Define $T^\circ$ to be the subset of $S^\circ$ on
which $v>v'$, and $T$ to be the closure of $T^\circ$. Suppose 
for a contradiction that $T$ is nonempty. Then $v>v'$ on $T^\circ$ 
and $v=v'$ on $\pd T$. As $v,v'$ are real analytic in $S^\circ$ by
assumption, it follows that $T$ has piecewise-smooth boundary.
Applying Proposition \ref{un8prop3} then shows that $v=v'$ on $T$,
a contradiction. Hence $T$ is empty, and $v\le v'$ on~$S$.
\end{proof}

If $v,v'\in C^2(S)$ satisfy \eq{un8eq1} then $Q(v)=Q(v')=0$ and
$v,v'$ are real analytic in $S^\circ$ by Proposition \ref{un8prop2}.
So we have:

\begin{cor} Let\/ $S$ be a domain in $\R^2$, let\/ $a\ne 0$,
and suppose $v,v'\in C^2(S)$ satisfy \eq{un8eq1} on $S$.
If\/ $v\le v'$ on $\pd S$ then $v\le v'$ on~$S$.
\label{un8cor}
\end{cor}

In particular, if $v\vert_{\pd S}=v'\vert_{\pd S}$ this gives
$v\le v'$ and $v'\le v$ on $S$, so that $v=v'$. This implies
uniqueness of solutions of the Dirichlet problem for \eq{un8eq1}.
Here is an analogue of Corollary \ref{un8cor} but with strict
inequalities, proved using a different method.

\begin{prop} Let\/ $S$ be a domain in $\R^2$, let\/ $a\ne 0$,
and suppose $v,v'\in C^2(S)$ satisfy \eq{un8eq1} on $S$.
If\/ $v<v'$ on $\pd S$ then $v<v'$ on~$S$.
\label{un8prop5}
\end{prop}

\begin{proof} Suppose for a contradiction that there exists 
$(b,c)\in S^\circ$ with $v(b,c)=v'(b,c)$. By Proposition
\ref{un8prop1} there exist $u,u'\in C^2(S)$, unique up to
addition of constants, such that $u,v$ and $u',v'$ satisfy 
\eq{un4eq5}. Choose the constants such that~$u(b,c)=u'(b,c)$.

Now apply Theorem \ref{un6thm2} to $(u,v)$ and $(u',v')$. As 
$v<v'$ on $\pd S$, the winding number of $(u,v)-(u',v')$
about 0 along $\pd S$ is zero, since $(u,v)-(u',v')$ is
confined to a half-plane and cannot go round $(0,0)$. But
$(u,v)-(u',v')$ has at least one zero in $S^\circ$, at $(b,c)$.
This is a contradiction. Therefore $v\ne v'$ in $S^\circ$, and 
by continuity and connectedness we have $v<v'$ on~$S$.
\end{proof}

\subsection{The Dirichlet problem for $v$}
\label{un82}

We now show that the Dirichlet problem for $v$ is uniquely
solvable in arbitrary domains $S$ in~$\R^2$.

\begin{thm} Let\/ $S$ be a domain in $\R^2$. Then whenever $a\ne 0$, 
$k\ge 0$, $\al\in(0,1)$ and\/ $\phi\in C^{k+2,\al}(\pd S)$ there 
exists a unique solution $v\in C^{k+2,\al}(S)$ of\/ \eq{un8eq1} 
with\/ $v\vert_{\pd S}=\phi$. Fix a basepoint\/ $(x_0,y_0)\in S$.
Then there exists a unique $u\in C^{k+2,\al}(S)$ with\/ $u(x_0,y_0)=0$
such that $u,v$ satisfy \eq{un4eq5}. Furthermore, $u,v$ are real
analytic in~$S^\circ$.
\label{un8thm1}
\end{thm}

\begin{proof} Observe that the operator $Q$ of \eq{un8eq1} is of
the form \eq{un3eq6}, with coefficients $a^{ij}$ depending on $y$
and $v$ but not on $\pd v$, and
\begin{equation*}
b\bigl((x,y),v,\pd v\bigr)=-\frac{v}{\bigl(v^2+y^2+a^2\bigr)^{3/2}}
\Bigl(\frac{\pd v}{\pd x}\Bigr)^2.
\end{equation*}
As $\bmd{v(v^2+y^2+a^2)^{-3/2}}\le a^{-2}$ the condition
$\bmd{b(x,u,p)}\le C\ms{p}$ in Theorem \ref{un3thm4} holds with
$C=a^{-2}$, and the condition $v\,b\bigl((x,y),v,p\bigr)\le 0$ 
for all $\bigl((x,y),v,p\bigr)\in S\t\R\t\R^2$ also clearly holds. 

Thus Theorem \ref{un3thm4} applies, and there exists $v$ in
$C^{k+2,\al}(S)$ satisying \eq{un8eq1} with $v\vert_{\pd S}=\phi$.
Corollary \ref{un8cor} shows that $v$ is unique. Using the condition
$u(x_0,y_0)=0$ to fix the additive constant, Proposition \ref{un8prop1}
shows that there exists a unique $u\in C^2(S)$ with $u(x_0,y_0)=0$
such that $u,v$ satisfy \eq{un4eq5}. But \eq{un4eq5} shows that
$\frac{\pd u}{\pd x},\frac{\pd u}{\pd y}\in C^{k+1,\al}(S)$ as
$v\in C^{k+2,\al}(S)$, so $u\in C^{k+2,\al}(S)$. Finally, Proposition
\ref{un4prop3} shows that $u,v$ are real analytic in~$S^\circ$.
\end{proof}

Combining the theorem with Proposition \ref{un4prop1}, we again get
an existence and uniqueness result for $\U(1)$-invariant SL 3-folds
in $\C^3$ satisfying certain boundary conditions, but different
boundary conditions to those in \S\ref{un73}. An analogue of Theorem
\ref{un8thm1} for the case $a=0$ will be given in \cite[Th.~6.1]{Joyc6},
which shows that \eq{un8eq1} has a unique {\it weak}\/ solution
$v\in C^0(S)$ with~$v\vert_{\pd S}=\phi$.

In Theorem \ref{un7thm2} we restricted $S$ to be a strictly convex
domain, but Theorem \ref{un8thm1} works for general domains. The
basic reason for this is that in the Dirichlet problem for $v$ we
automatically get an a priori estimate for $\cnm{v}{0}$, which
implies positive upper and lower a priori bounds
for~$(v^2+y^2+a^2)^{-1/2}$.

Hence, in the Dirichlet problem for $v$ we know in advance that 
\eq{un8eq1} is {\it uniformly elliptic}. However, in the Dirichlet 
problem for $f$ we need an a priori bound for $\cnm{\frac{\pd f}{\pd 
x}}{0}$ to make \eq{un7eq1} uniformly elliptic, and we assume $S$ is
strictly convex to prove such a bound.

By analogy with Theorem \ref{un7thm3}, we can also prove:

\begin{thm} Let\/ $S$ be a domain in $\R^2$, $k\ge 0$, $\al\in(0,1)$
and\/ $(x_0,y_0)\in S$. Then the map $C^{k+2,\al}(\pd D)\t\bigl(\R\sm
\{0\}\bigr)\ra C^{k+2,\al}(D)^2$ taking $(\phi,a)\mapsto(u,v)$ is
continuous, where $(u,v)$ is the unique solution of\/ \eq{un4eq5}
with\/ $v\vert_{\pd D}=\phi$ and\/ $u(x_0,y_0)=0$ constructed in
Theorem~\ref{un8thm1}.
\label{un8thm2}
\end{thm}

Presumably the map $(\phi,a)\mapsto(u,v)$ is also smooth. An
extension of Theorem \ref{un8thm2} to include the case $a=0$ is
given in \cite[Th.~6.2]{Joyc6}, but with the $C^0$ rather than
the $C^{k+2,\al}$ topology on~$u,v$.

\subsection{Winding number results for $v$}
\label{un83}

As in \S\ref{un74}, we will now extend some of the `winding 
number' results of \S\ref{un6} to the situation of this 
section. Here is the analogue of Morse function for~$v$.

\begin{dfn} Let $S$ be a domain in $\R^2$, and let $v\in C^1(\pd S)$.
Choose a smooth, positively oriented parametrization $\th$ for
$\pd S$ as in \eq{un7eq6}, and regard $v$ as a function of $\th$.
We call $v$ {\it transverse} if $v=0$ at only finitely many points
in $\pd S$, and $\frac{\d v}{\d\th}\ne 0$ at each of these points.

This definition is independent of parametrization $\th$, and
transverse functions are an open dense subset of $C^1(\pd S)$.
Also, each zero of $v$ is either {\it increasing}, with 
$\frac{\d v}{\d\th}>0$, or {\it decreasing}, with
$\frac{\d v}{\d\th}<0$, and there are the same number of 
each, so $f$ has exactly $l$ increasing and $l$ decreasing
zeroes for some~$l\ge 0$.
\label{un8def1}
\end{dfn}

Here is the analogue of Theorem \ref{un7thm4}, with a similar proof.

\begin{thm} Let\/ $S$ be a domain in $\R^2$, let\/ $a\ne 0$, and 
let\/ $u_1,v_1\in C^1(S)$ and\/ $u_2,v_2\in C^1(S)$ be solutions 
of\/ \eq{un4eq5}. Suppose that\/ $(v_1-v_2)\vert_{\pd S}$ is 
transverse with\/ $2l$ zeroes. Then $(u_1,v_1)-(u_2,v_2)$ has $n$ 
zeroes in $S^\circ$ with multiplicities $k_1,\ldots,k_n$ and\/ $m$ 
zeroes on $\pd S$, where $n,m\ge 0$ and\/ $k_i\ge 1$, 
and\/~$\sum_{i=1}^nk_i+m\le l$.
\label{un8thm3}
\end{thm}

\begin{proof} Suppose $(u_1,v_1)=(u_2,v_2)$ at $(x_0,y_0)$ in $\pd S$. 
Then $v_1-v_2=0$ at $(x_0,y_0)$, so $(x_0,y_0)$ must be one of the $2l$ 
zeroes of $v_1-v_2$ on $\pd S$. Thus $m$ is finite. Let $m_1\ge 0$ of 
the $m$ zeroes of $(u_1,v_1)-(u_2,v_2)$ on $\pd S$ be {\it increasing} 
zeroes of $v_1-v_2$ on $\pd S$, and $m_2\ge 0$ be {\it decreasing} 
zeroes, where $m_1+m_2=m$. As in the proof of Theorem \ref{un7thm4} we
find that $(u_1,v_1)-(u_2,v_2)$ has finitely many zeroes in $S^\circ$, 
so let there be $n$ zeroes, with multiplicities~$k_1,\ldots,k_n$.

Let $\ep>0$ be small. Then $(u_1+\ep,v_1)$ also satisfies \eq{un4eq5},
so we can consider the zeroes of $(u_1+\ep,v_1)-(u_2,v_2)$ in $S$.
One can use the ideas of \S\ref{un6} to show that close to the $n$ 
zeroes of $(u_1,v_1)-(u_2,v_2)$ in $S^\circ$ there will be $n'\ge n$ 
zeroes of $(u_1+\ep,v_1)-(u_2,v_2)$ with multiplicities $k_1',\ldots,
k_{n'}'$, where~$\sum_{i=1}^{n'}k_i'=\sum_{i=1}^nk_i$.

If $(x_0,y_0)$ is one of the $m_1$ increasing zeroes of $v_1-v_2$ 
on $\pd S$, as $\ep>0$ is small one can show that $(u_1+\ep,v_1)-
(u_2,v_2)$ has a zero in $S^\circ$ near $(x_0,y_0)$. If $(x_0,y_0)$ 
is one of the $m_2$ decreasing zeroes of $v_1-v_2$ on $\pd S$, there 
is no zero of $(u_1+\ep,v_1)-(u_2,v_2)$ in $S^\circ$ near $(x_0,y_0)$.
So, by Theorem \ref{un6thm2} the winding number of $(u_1+\ep,v_1)
-(u_2,v_2)$ about 0 along $\pd S$ is~$k'=\sum_{i=1}^{n'}k_i'+m_1$.

Now $(u_1+\ep,v_1)-(u_2,v_2)$ crosses the $x$-axis exactly $2l$
times on $\pd S$, at the zeroes of $v_1-v_2$. However, at the
$m_2$ decreasing zeroes $(u_1+\ep,v_1)-(u_2,v_2)$ crosses the
$x$-axis at $(\ep,0)$ in the negative sense winding round 0.
So it is not difficult to see that the winding number $k'$
satisfies $k'\le l-m_2$. The theorem then follows from the
equations $k'=\sum_{i=1}^{n'}k_i'+m_1$, $\sum_{i=1}^{n'}k_i'=
\sum_{i=1}^nk_i$ and~$m_1+m_2=m$.
\end{proof}

The theorem can be used in conjunction with Theorem \ref{un8thm1},
the solution of the Dirichlet problem for $v$. In this case, we 
would know $v_1,v_2$ on $\pd S$ explicitly, but would otherwise 
know little about the $u_j$ or $v_j$. The theorem tells us something 
about $(u_1,v_1)$ and $(u_2,v_2)$, using only the known boundary 
values of~$v_1,v_2$.

\subsection{A maximum principle for $\frac{\pd v}{\pd x}$}
\label{un84}

Finally we show that if $v$ satisfies \eq{un8eq1} on $S$ then
$\bmd{\frac{\pd v}{\pd x}}$ is maximum on~$\pd S$.

\begin{prop} Let\/ $S$ be a domain in $\R^2$, let\/ $a\ne 0$,
and suppose $v\in C^2(S)$ satisfies \eq{un8eq1}. Then the maximum
of\/  $\bmd{\frac{\pd v}{\pd x}}$ is achieved on~$\pd S$.
\label{un8prop6}
\end{prop}

\begin{proof} As $v$ satisfies \eq{un8eq1} it is real analytic in
$S^\circ$ by Proposition \ref{un8prop2}, and satisfies \eq{un8eq2}.
Taking the derivative $\frac{\pd}{\pd x}$ of \eq{un8eq2} in
$S^\circ$ and rearranging gives
\begin{align*}
L\Bigl(\frac{\pd v}{\pd x}\Bigr)&=
\bigl(v^2+y^2+a^2\bigr)^{-3/2}\Bigl(\frac{\pd v}{\pd x}\Bigr)^3,
\quad\text{where}\\
L(g)&=\bigl(v^2+y^2+a^2\bigr)^{-1/2}\frac{\pd^2g}{\pd x^2}
+2\frac{\pd^2g}{\pd y^2}-3\bigl(v^2+y^2+a^2\bigr)^{-3/2}v
\frac{\pd v}{\pd x}\cdot\frac{\pd g}{\pd x}.
\end{align*}
Then $L$ is a linear elliptic operator of the form \eq{un3eq1},
with~$c(x)\equiv 0$.

Suppose that $\frac{\pd v}{\pd x}$ has a positive maximum achieved
at $(x_0,y_0)\in S^\circ$, with $\frac{\pd v}{\pd x}(x_0,y_0)=M>0$
say, and $\frac{\pd v}{\pd x}<M$ on $\pd S$. Let $\ep\in(0,M)$ be
generic and small enough that $\frac{\pd v}{\pd x}<M-\ep$ on $\pd S$,
and define $T=\bigl\{(x,y)\in S:\frac{\pd v}{\pd x}\ge M-\ep\bigr\}$.
Then as $\ep$ is generic $T$ lies in $S^\circ$ and is compact with
smooth boundary, and $\frac{\pd v}{\pd x}=M-\ep$ on $\pd T$. Also
$L\bigl(\frac{\pd v}{\pd x}\bigr)>0$ on $T$, as $\frac{\pd v}{\pd x}
\ge M-\ep>0$ on~$T$.

Applying Theorem \ref{un3thm2} shows that the maximum of
$\frac{\pd v}{\pd x}$ on $T$ is achieved on $\pd T$. This
contradicts $(x_0,y_0)\in T^\circ$, $\frac{\pd v}{\pd x}(x_0,y_0)
=M$ and $\frac{\pd v}{\pd x}=M-\ep$ on $\pd T$. Thus, if
$\frac{\pd v}{\pd x}$ has a positive maximum it is achieved on
$\pd S$. Similarly, if $\frac{\pd v}{\pd x}$ has a negative
minimum it is achieved on $\pd S$. Thus the maximum of
$\bmd{\frac{\pd v}{\pd x}}$ is achieved on~$\pd S$.
\end{proof}

\end{document}